\documentclass[reqno]{amsart}
\usepackage[usenames,dvipsnames]{pstricks}
\usepackage{epsfig}
\usepackage[english]{babel}
\usepackage[latin1]{inputenc}
\usepackage{graphicx}
\usepackage{subfigure}
\usepackage{xspace} % To get the right spacings in front of : and so on
\usepackage{amsmath,amssymb}
\usepackage{fancyhdr}

\makeatletter
\newcommand{\sech}{\mathop{\operator@font sech}}
\newcommand{\sign}{\mathop{\operator@font sign}}
\makeatother

\numberwithin{equation}{section}

\begin{document}
%\begin{frontmatter}
\title{Notes on the numerical solution of the Benjamin equation}

\author[V. A. Dougalis]{Vassilios A. Dougalis}
\address{Mathematics Department, University of Athens, 15784
Zographou, Greece \and Institute of Applied \& Computational
Mathematics, FO.R.T.H., 71110 Heraklion, Greece}
\email{doug@math.uoa.gr}

\author[A. Duran]{Angel Duran}
\address{ Applied Mathematics Department,  University of
Valladolid, 47011 Valladolid, Spain}
\email{angel@mac.uva.es}

\author[D. Mitsotakis]{Dimitrios Mitsotakis}
\address{Applied Mathematics Department, University of California, Merced, 5200 North Lake rd. Merced, CA 95343 , USA}
\email{dmitsot@gmail.com}
\urladdr{http://dmitsot.googlepages.com/}

\subjclass[2010]{76B15 (primary), 65M60, 65M70 (secondary)}

\keywords{Benjamin equation, Solitary waves, Hybrid Finite Element-Spectral method}

\begin{abstract}
In this paper we consider the  Benjamin equation, a partial differential equation  that models  one-way propagation of  long internal waves of small amplitude along the interface of two fluid layers under the effects of gravity and surface tension. We solve the periodic initial-value problem for the Benjamin equation numerically by a new  fully discrete hybrid finite-element / spectral scheme, which we first validate by pinning down its accuracy and stability properties.  After testing the evolution properties of the scheme in a study of propagation of single - and multi-pulse solitary waves of the Benjamin equation, we use it in an exploratory mode to illuminate phenomena such as overtaking collisions of solitary waves, and the stability of single-, multi-pulse  and `depression' solitary waves.

\end{abstract}

\maketitle

\section{Introduction}
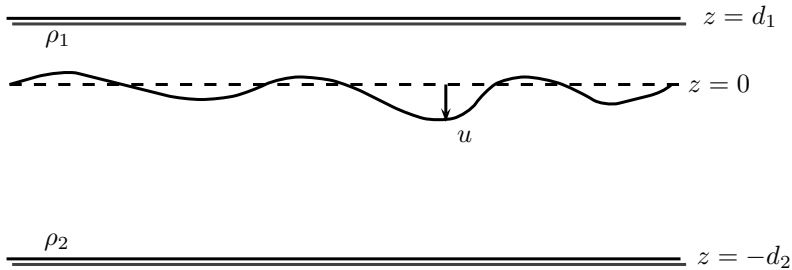
\begin{figure}[h]\label{F1a}
\centering
\scalebox{1} % Change this value to rescale the drawing.
{
\begin{pspicture}(0,-1.81)(10.54,1.81)
\psline[linewidth=0.04cm,shadow=true](0.0,1.59)(8.96,1.59)
\psline[linewidth=0.04cm,shadow=true](0.0,-1.61)(8.96,-1.61)
\psline[linewidth=0.04cm,linestyle=dashed,dash=0.16cm 0.16cm](0.04,0.71)(8.94,0.71)
\pscustom[linewidth=0.04]
{
\newpath
\moveto(0.04,0.71)
\lineto(0.39,0.81)
\curveto(0.565,0.86)(0.84,0.885)(0.94,0.86)
\curveto(1.04,0.835)(1.24,0.785)(1.34,0.76)
\curveto(1.44,0.735)(1.64,0.685)(1.74,0.66)
\curveto(1.84,0.635)(2.04,0.585)(2.14,0.56)
\curveto(2.24,0.535)(2.465,0.51)(2.59,0.51)
\curveto(2.715,0.51)(2.94,0.535)(3.04,0.56)
\curveto(3.14,0.585)(3.34,0.66)(3.44,0.71)
\curveto(3.54,0.76)(3.765,0.81)(3.89,0.81)
\curveto(4.015,0.81)(4.24,0.785)(4.34,0.76)
\curveto(4.44,0.735)(4.64,0.66)(4.74,0.61)
\curveto(4.84,0.56)(5.04,0.46)(5.14,0.41)
\curveto(5.24,0.36)(5.44,0.285)(5.54,0.26)
\curveto(5.64,0.235)(5.84,0.235)(5.94,0.26)
\curveto(6.04,0.285)(6.19,0.385)(6.24,0.46)
\curveto(6.29,0.535)(6.415,0.66)(6.49,0.71)
\curveto(6.565,0.76)(6.765,0.81)(6.89,0.81)
\curveto(7.015,0.81)(7.265,0.76)(7.39,0.71)
\curveto(7.515,0.66)(7.715,0.56)(7.79,0.51)
\curveto(7.865,0.46)(8.04,0.435)(8.14,0.46)
\curveto(8.24,0.485)(8.44,0.535)(8.54,0.56)
\curveto(8.64,0.585)(8.765,0.635)(8.84,0.71)
}
\rput(9.47,0.715){$z=0$}
\usefont{T1}{ptm}{m}{n}
\rput(9.75,1.615){$z=d_1$}
\usefont{T1}{ptm}{m}{n}
\rput(9.81,-1.585){$z=-d_2$}
\psline[linewidth=0.04cm,arrowsize=0.05291667cm 2.0,arrowlength=1.4,arrowinset=0.4]{->}(5.84,0.71)(5.84,0.21)
\usefont{T1}{ptm}{m}{n}
\rput(6.09,0.015){$u$}
\usefont{T1}{ptm}{m}{n}
\rput(0.68,1.315){$\rho_1$}
\usefont{T1}{ptm}{m}{n}
\rput(0.67,-1.385){$\rho_2$}
\end{pspicture}
}
\caption{Interfacial gravity-capillary waves}
\end{figure}

In this paper we will consider the {\it Benjamin equation}
\begin{equation}\label{E11}
u_t+\alpha u_x+\beta u u_x-\gamma \mathcal{H}u_{xx}-\delta u_{xxx}=0,
\end{equation}
where $u=u(x,t), x\in \mathbb{R}, t\geq 0, \alpha,\beta,\gamma,\delta$ are positive constants, and $\mathcal{H}$ denotes the Hilbert transform defined on the real line as
$$\mathcal{H}f(x):=\frac{1}{\pi}p.v.\int_{-\infty}^{\infty}\frac{f(y)}{x-y}\,dy$$ or through its Fourier transform as $$\widehat{\mathcal{H}f}(k)=-{\rm i}\sign(k)\widehat{f}(k), k\in\mathbb{R}.$$
The Benjamin equation, cf. \cite{B1,B2,ABR}, is a model for {\it internal} waves propagating under the effect of gravity and surface tension in the positive $x$-direction along the interface of a two-dimensional system of two homogeneous layers of incompressible, inviscid fluids consisting at rest of a thin layer of fluid 1 of depth $d_{1}$ and density $\rho_{1}$ lying above a layer of fluid 2 of very large depth $d_{2}\gg d_{1}$ and density $\rho_{2}>\rho_{1}$. The upper layer is bounded above by a horizontal \lq rigid lid\rq\, and the lower layer is bounded below by an impermeable horizontal bottom, as in Figure \ref{F1a}.

It is further assumed that the following physical regime of interest is to be modelled: Let $a$ be a typical amplitude and $\lambda$ a typical wavelength of the interfacial wave. The parameters $\epsilon=a/d_{1}$ and $\mu=d_{1}^{2}/\lambda^{2}$ are assumed to be small and satisfy $\mu\sim\epsilon^{2}\ll 1$; it is also assumed that capillarity effects along the interface are not negligible. Under these assumptions (\ref{E11}) was derived in \cite{B1} from the two-dimensional, two-layer Euler equations in the presence of interface surface tension by dispersion relation arguments. The variables in (\ref{E11}) are nondimensional and scaled, and the coefficients are given by $$\alpha=\sqrt{\frac{\rho_{2}-\rho_{1}}{\rho_{1}}}, \quad \beta=\frac{3}{2}\alpha \epsilon, \quad\gamma=\frac{1}{2}\alpha\sqrt{\mu}\frac{\rho_{2}}{\rho_{1}}, \quad\delta=\frac{\alpha T}{2g\lambda^{2}(\rho_{2}-\rho_{1})},$$ where $T$ is the interfacial surface tension and $g$ the acceleration of gravity. The variables $x$ and $t$ are proportional to distance along the channel and time, respectively, and $u(x,t)$ denotes the downward vertical displacement of the interface from its level of rest at $(x,t)$. The interfacial surface tension $T$ is assumed to be much larger than $g(\rho_{2}-\rho_{1})d_{1}^{2}$. (For a further discussion of the physical regime of validity of (\ref{E11}) cf. \cite{ABR}.) Note that if the parameter $\delta$ is taken equal to zero, (\ref{E11}) reduces to the Benjamin-Ono (BO) equation, \cite{B0,O}, while, if we put $\gamma=0$ we obtain the KdV equation with negative dispersion coefficient.

It is well known, cf. \cite{B1}, that sufficiently smooth solutions of (\ref{E11}) that vanish suitably at infinity preserve the functionals
\begin{eqnarray}
m(u)&=&\int_{-\infty}^{\infty} udx,\label{E12}\\
I(u)&=&\frac{1}{2}\int_{-\infty}^{\infty} u^{2}dx,\label{E13}\\
E(u)&=&\int_{-\infty}^{\infty}
\left(\frac{\beta}{6}u^{3}-\frac{1}{2}\gamma
u\mathcal{H}u_{x}+\frac{1}{2}\delta u_{x}^{2}\right)dx.\label{E14}
\end{eqnarray}
Global well-posedness in $L^{2}$ for the Cauchy problem and also for the periodic initial-value problem for (\ref{E11}) was established in \cite{L}.

In this paper we will study (\ref{E11}) numerically, paying particular attention to properties of its {\it solitary-wave} solutions. These are travelling-wave solutions of the form $u(x,t)=\varphi(x-c_{s}t), c_{s}>0$, such that $\varphi$ and its derivatives tend to zero  as $\xi=x-c_{s}t$ approaches $\pm\infty$. Substituting this expression in (\ref{E11}) and integrating once we obtain
\begin{eqnarray}
(\alpha-c_{s})\varphi
+\frac{\beta}{2}\varphi^{2}-\gamma \mathrm{H}\varphi-\delta
\varphi^{\prime\prime}=0,\label{E15}
\end{eqnarray}
where ${}^{\prime}=d/d\xi$, and the operator $\mathrm{H}$ is defined by $\mathrm{H}:=\mathcal{H}\partial_{x}$, i.~e. by $\widehat{\mathrm{H}f}(k)=|k|\widehat{f}(k), k\in \mathbb{R}$. We will assume that $\alpha-c_{s}>0$.

If we perform the change of variables
\begin{eqnarray*}
\varphi(\xi)=-\frac{2(\alpha-c_{s})}{\beta}\psi(z),\quad z=\sqrt{\frac{\alpha-c_{s}}{\delta}}\xi,
\end{eqnarray*}
in (\ref{E15}), we see that the solitary-wave profile $\psi(z)$ satisfies the ordinary differential equation (ode)
\begin{eqnarray}
\psi-2\tilde{\gamma}\mathrm{H}\psi-\psi_{zz}-\psi^{2}=0,\quad z\in\mathbb{R},\label{E16}
\end{eqnarray}
where
\begin{eqnarray}
\tilde{\gamma}=\frac{\gamma}{2\sqrt{\delta(\alpha-c_{s})}}.\label{E17}
\end{eqnarray}
This change of variables and the resulting equation (\ref{E16}) was used in \cite{B1,B2}, and \cite{ABR}. (In these references $\tilde{\gamma}$ is denoted by $\gamma$.) In his papers Benjamin showed, using degree theory, that for each $\tilde{\gamma}\in [0,1)$, there exists a solution $\psi$ of (\ref{E16}) which is an even function of $z$ with $\psi(0)=\max_{z\in\mathbb{R}}\psi(z)>0$. He also argued by formal asymptotics that for each $\tilde{\gamma}\in [0,1)$ there is a bounded interval centered at $z=0$, in which $\psi$ oscillates (with the number of oscillations increasing as $\tilde{\gamma}$ approaches $1$), while outside this interval he concluded in \cite{B2} that $|\psi|$ decays like $1/z^{2}$. In addition, in the same paper he outlined an orbital stability theory for these solitary waves for small $\tilde{\gamma}$. In \cite{ABR} a complete theory of existence and orbital stability of the solitary waves for small $\tilde{\gamma}$ was presented, based on the implicit function theorem, perturbation theory of operators, and the fact that $\tilde{\gamma}=0$ corresponds to solitary waves of the KdV equation. Further issues of existence and rigorous asymptotics of the solitary waves of (\ref{E11}) and related equations were explored in \cite{CB}. In \cite{A} concentration compactness arguments were used to establish existence and a weaker version of stability of the solitary waves of (\ref{E11}) for $0<\tilde{\gamma}<1$.

In this paper we will employ the solitary-wave equation in the form (\ref{E15}). As a result, normally the solitary waves will have negative maximum excursions from their level of rest.

Since explicit formulas for the solitary waves of the Benjamin equation are not known (except when one of $\gamma$ or $\delta$ is set equal to zero), one must resort to approximate techniques for their construction. The presence of the nonlocal terms in (\ref{E11}) and (\ref{E15}), which have a handy Fourier representation in the periodic case as well, naturally suggests using spectral-type methods for approximating their solutions. The preceding discussion of the Benjamin equation applies to its associated Cauchy problem on $\mathbb{R}$. Solving it numerically requires posing it on a finite $x$-interval $[-L,L]$ with, say, periodic boundary conditions, assuming $2L$-periodic initial data. In case solitary waves, their generation and interactions, are the focus of interest, one should take into account that they decay quadratically. Consequently, the interval $[-L,L]$ should be taken sufficiently large in some experiments to ensure that the numerical solution in the temporal range of interest remains sufficiently small at the endpoints so that the simulations give valid approximations of the solutions of the Cauchy problem.

In \cite{ABR} the equation (\ref{E16}) was discretized in space by a pseudospectral technique and the resulting nonlinear system of equations for the Fourier coefficients of $\psi=\psi_{\tilde{\gamma}}$ for a desired value of $\tilde{\gamma}\in (0,1)$ was solved by an incremental continuation method. This entailed defining a homotopic path $\tilde{\gamma}_{0}=0<\tilde{\gamma_{1}}<\ldots<\tilde{\gamma}_{M}=\tilde{\gamma}$, starting from the known profile of a solitary wave $\psi_{\tilde{\gamma}_{0}}$ of the KdV equation with a given speed $c_{s}$, and computing $\psi_{\tilde{\gamma}_{j+1}}$, given $\psi_{\tilde{\gamma}_{j}}$, by Newton's method. With this technique the authors of \cite{ABR} were able to construct approximate solutions of (\ref{E16}) that were even functions with a positive absolute maximum at $z=0$. As $\tilde{\gamma}$ approached $1$ the oscillating tails of the solitary wave became more prominent and the maximum value of the wave decreased. It was found that the length of the intervals between consecutive zeros of the oscillating tails was quite close to the value predicted by the asymptotic analysis of \cite{B2}.

In \cite{KB} the authors solved numerically the periodic initial-value problem for the Benjamin equation using a pseudospectral (collocation) method in space coupled with a second-order time-stepping procedure. They confirmed that {\it resolution} of suitable general initial profiles into a number of solitary waves plus a dispersive tail (a phenomenon that has been observed in other nonlinear dispersive wave equations) also occurs in the case of the Benjamin equation. They specifically studied the resolution of initial Gaussian profiles into solitary waves contrasting it with the analogous resolution observed in the case of two BO-type equations. In some cases they observed, in addition to detached solitary waves, the emergence of clusters (pairs, triplets, etc.) of \lq orbiting\rq\ solitary waves that interacted among themselves. They conjectured that these structures would eventually separate into distinct solitary waves. They also constructed approximate solitary waves, using the resolution property, by truncating and iteratively \lq cleaning\rq\ a separated solitary wave as has been frequently done in numerical studies of other nonlinear dispersive wave equations. (Of course in this manner one does not have in general {\it a priori} knowledge of the speed $c_{s}$ or the value of $\tilde{\gamma}$ of the emerging solitary wave.) They used two such approximate  solitary waves of different speeds to study their overtaking collision and observed that the interaction was not elastic, a fact indicating that the Benjamin equation is not integrable.

In \cite{CA}, the authors considered solitary waves of the Benjamin equation and compared them to solitary waves of the full Euler equations for interfacial flows in the presence of surface tension when the parameters of the problem are close to the Benjamin equation regime of validity and also farther from it. The numerical scheme they used for approximating solitary waves of the Benjamin equation was based on a hybrid spatial discretization that employed fourth-order finite differences on a uniform grid for the derivatives, and the discrete Fourier transform for the nonlocal term. The resulting nonlinear system of equations was solved again by a continuation-Newton technique. The temporal discretization of the periodic initial-value problem for the Benjamin equation was effected by an explicit predictor-corrector scheme. They identified another branch of solitary wave solutions of the Benjamin equation, the \lq depression\rq\ solitary waves (resembling analogous solutions of the Euler equations), and tested their stability by using them as initial values in their fully discrete scheme for the time-dependent equation. They observed that the initial profile propagated without change for some time, gradually developed an instability due to the perturbative effect of the numerical scheme, and resolved itself into two pulses resembling usual (\lq elevation\rq) solitary waves of the Benjamin equation plus small-amplitude dispersive oscillations. (A linearized stability analysis, also performed in \cite{CA}, yields that the depression solitary waves are linearly unstable.)

In a recent paper \cite{DDM1}, we made a study of several incremental continuation techniques for approximating solitary waves of the Benjamin equation that satisfy (\ref{E15}). (The values of $\alpha, \beta, \delta$ and $c_{s}$ were fixed, and $\gamma$ was used as continuation parameter.) A standard pseudospectral (collocation) method yielded the underlying discrete nonlinear system. We found that Newton's method, combined with a suitably preconditioned conjugate gradient technique for solving the attendant linear system at each Newton iteration, was the generally most efficient technique of implementing the incremental step and produced very accurate approximations of the solitary waves for $0\leq \gamma <1$. With this method we also computed other branches of solutions of (\ref{E15}), namely {\it multi-pulse} solitary waves, by starting the homotopy path from linear combinations of solitary waves of the KdV equation. We verified the accuracy of these profiles as travelling waves of the Benjamin equation by using them as initial values in a full discretization of the periodic initial-value problem for (\ref{E11}) and integrating forward in time. The solver combined the pseudospectral spatial discretization with the third-order accurate two-stage DIRK time-stepping technique, modified to preserve discrete analogs of the invariants (\ref{E12}) and (\ref{E13}). It was found that several quantities of interest, such as the speed, the amplitude and the third invariant (\ref{E14}) of the discrete travelling waves, were preserved to very high accuracy, lending confidence in the validity of this technique for computing solitary waves.

In the paper at hand we continue our numerical study of the Benjamin equation. We construct and test numerically a new, efficient time-stepping method based on a spectral-finite element {hybrid} spatial discretization combined with a fourth-order implicit Runge-Kutta scheme for time-stepping. This method is used to explore properties of solitary-wave solutions of (\ref{E11}), such as their generation, interaction and stability.

Much of numerical work with spectral-type methods for one-dimensional, nonlocal, nonlinear dispersive wave equations has been centered around the Benjamin-Ono (BO), \cite{B0,O}, and the Intermediate Long Wave (ILW) equation, \cite{J,ABFS}. Early computational work was reviewed in \cite{PD1}; here we mention only the rigorous convergence results known to us. In \cite{PD2} $L^{2}-$error estimates were derived for the standard Fourier-Galerkin semidiscretization of the BO and ILW equations. If the number of Fourier modes is $2N+1$ and the initial value is $2L-$periodic and belongs to the periodic Sobolev space $H_{p}^{r}$, the $L^{2}$-error bounds derived in \cite{PD2} are of $O(N^{1-r})$. In addition, the full discretization of the semidiscrete system of ode's with the explicit leap-frog scheme is shown in \cite{PD2} to have an $L^{2}-$error bound of $O(N^{1-r}+\Delta t^{2})$ under the stability restriction that $N^{2}\Delta t\leq C$ for a sufficiently small constant $C$; here $\Delta t$ is the time step. For a class of equations with the same nonlocal terms and more general nonlinear terms  it was subsequently shown in \cite{DM1} that the error of the Fourier-Galerkin semidiscretization is of optimal order $O(N^{1/2-r})$ in $H^{1/2}_{p}$. In the same paper the semidiscrete problem was discretized in time in the manner suggested in \cite{CK}, i.~e. using as a basis the leap-frog method coupled with implicit Crank-Nicolson differencing of the linear dispersive term. This explicit-implicit time-stepping scheme may be implemented efficiently in Fourier space and does not require solving linear systems of equations; as shown in \cite{DM1} it has an error bound of $O(N^{1/2-r}+\Delta t^{2})$ in $H^{1/2}_{p}$ under the mild stability condition $N^{1/2}\Delta t\leq C$ for some sufficiently small constant $C$. In addition, in \cite{DM2} the authors analyze the more efficient spectral collocation method (that was used in actual computations in \cite{PD1} and elsewhere,) for the BO and ILW equations, and prove that the associated semidiscrete problem converges with an $H^{1/2}_{p}-$error bound of $O(N^{3/2-r})$.

A different type of method for the BO equation was constructed and analyzed in \cite{TM}. It consists of a Crank-Nicolson time-stepping scheme that is coupled with a spatial discretization in which the nonlinear term is approximated by conservative differencing and the nonlocal term is discretized in physical space by the midpoint quadrature formula, which is then interpreted as a discrete convolution and computed by the discrete Fourier transform. Since the fully discrete scheme is implicit, a nonlinear system of equations has to be solved at each time step. This system is linearized by a simple iterative scheme in which the nonlinear term is lagged backwards in time and the linear part is trivial to invert in Fourier space, as e.~g. in \cite{CK}. The overall method is shown to be of second-order accuracy in $L^{2}$ in space and time.

In the present paper the numerical scheme that we use is a {hybrid} finite element-spectral method. We consider the periodic initial-value problem for (\ref{E11}) and discretize it in space by the Galerkin method using smooth periodic splines of order $r\geq 3$ on a uniform mesh with meshlength $h$. (Cubic splines, i.~e. $r=4$, are mainly used in the computations.) The nonlocal term is computed using a spectral approximation as described in Section \ref{sec2}. Then, the system of ode's representing the semidiscrete problem is discretized in time; we use as a base time-stepping scheme the two-stage, fourth-order accurate, Gauss-Legendre implicit Runge-Kutta method. This scheme has high accuracy and good stability properties and has previously been extensively used for the temporal discretization of stiff partial differential equations with a KdV term, cf. e.~g. \cite{BDKMc} and its references. We describe in detail the implementation of this fully discrete hybrid method and make a computational study of its accuracy and stability properties when it is applied to the Benjamin and Benjamin-Ono (i.~e. when $\delta$ is set to zero) equations. In addition, we validate the hybrid scheme by making a detailed comparison of the solutions that it produces with those of a standard fully discrete pseudospectral scheme in the case of three numerical experiments involving the propagation of solitary waves of the Benjamin and Benjamin-Ono equations.

In Section \ref{sec3} we review the continuation-conjugate gradient-Newton technique of \cite{DDM1} for generating single and multi-pulse solitary-wave solutions (i.~e. solutions of (\ref{E15})) of the Benjamin equation for various values of $\gamma$ with particular attention to values close to $1$. We use these numerical profiles as initial conditions in numerical evolution experiments with the hybrid scheme and investigate with various metrics their accuracy as travelling wave solutions of the Benjamin equation. Our conclusion from the numerical experiments of Sections \ref{sec2} and \ref{sec3} is that the hybrid scheme yields very accurate and stable approximations of solutions of the Benjamin equation, and in particular of the solitary waves for values of $\gamma\in (0,1)$ that can be taken quite close to $1$.

In Section \ref{sec4} we make a detailed computational study of overtaking (`one-way') collisions of solitary waves of the Benjamin equation and compare the inelastic character of these interactions with the analogous, \lq clean\rq\ interactions in the case of the integrable BO equation. Finally, in Section \ref{sec5} we explore issues of stability and instability of single-and multi-pulse solitary waves of the Benjamin equation under small and large perturbations. Our computational study confirms the stability of the single-pulse solitary waves for small and moderate values of $\gamma$ but is inconclusive for cases of $\gamma$ very close to $1$. The multi-pulse waves appear to be unstable and our experiments suggest that after an initial \lq orbiting\rq\ or \lq dancing\rq\ phase, they produce separated solitary waves. This confirms the conjecture of \cite{KB} that was mentioned previously. Finally, we examine the stability of the \lq depression\rq\ solitary waves and confirm the results of \cite{CA} regarding their instability.

In summary, the main contributions of the paper at hand are
\begin{itemize}
\item The construction of a novel, highly accurate, stable and efficient hybrid scheme that combines the accuracy of the spectral approximation of the nonlocal term and the accuracy of the spline discretization of the rest of the terms of the Benjamin equation with an accurate, unconditionally stable time stepping procedure which is effective in approximating highly stiff problems such as semidiscretizations of the Benjamin equation in the presence of the KdV term.
\item The validation of the accuracy of the numerically generated single- and multi- pulse solitary wave solutions by showing that when used as initial values of the hybrid scheme they produce highly accurate approximations to travelling wave solutions of the evolution problem. These approximate solitary waves were computed by a Fourier spatial discretization of the solitary wave ode (\ref{E15}) coupled with a continuation conjugate gradient-Newton nonlinear system solver that was proposed by the authors in \cite{DDM1} and can produce accurate solitary waves for any desired values of the speed $c_{s}>\alpha$ and $\gamma\in [0,1)$, avoiding the drawbacks of the iterative \lq cleaning\rq\ .
\item The illumination, by computational means, of important phenomena associated with solitary waves of nonlinear dispersive wave equations, such as their one-way interaction (overtaking collision) and stability properties in the case of the Benjamin equation.
\end{itemize}

In the paper, we denote , for integer $r\geq 0$, by $C_{p}^{r}$ the periodic functions, on $[-L,L]$ or $[0,2\pi]$ as the case may be, that belong to $C^{r}$. The inner product for real or complex-valued functions in $L^{2}$ is denoted by $(\cdot,\cdot)$ and the associated norm by $||\cdot ||$.
\section{The hybrid spectral-finite element scheme}
\label{sec2}
We consider the periodic initial-value problem for the Benjamin equation, i.~e. for $t\geq 0$ we seek a $2L-$periodic real function $u=u(x,t)$ such that
\begin{eqnarray}
& u_t+\alpha u_x+\beta u u_x-\gamma \mathcal{G}u_{xx}-\delta u_{xxx}=0,\quad x\in[-L,L], \quad t>0,  \label{E21}\\
& u(x,0)=u_0(x),\quad x\in [-L,L]\nonumber
\end{eqnarray}
where $u_{0}$ is a given smooth $2L-$periodic function and $\alpha,\beta, \gamma, \delta$ positive constants. The operator $\mathcal{G}$ is the Hilbert transform acting on $2L-$periodic functions; for the purposes of this section it will be represented by its principal-value integral form \cite{ABFS}
\begin{eqnarray}
\mathcal{G}f(x):=\frac{1}{2L}p.v.\int_{-L}^{L}{\it cot}\left(\frac{\pi (x-y)}{2L}\right)f(y)dy,\label{E22}
\end{eqnarray}
where $f$ is $2L-$periodic. In the sequel we will assume that the solution of (\ref{E21}) is sufficiently smooth. For simplicity, we assume that the problem (\ref{E21}) has been transformed onto the spatial interval $[0,2\pi]$.
\subsection{The semidiscrete hybrid scheme}
\label{sec21}
For integer $r\geq 3$ and an even integer $N$, let $h=2\pi/N$, $x_{j}=jh, j=0,\ldots,N$, and consider the finite dimensional spaces
$$S_N=\mbox{span}\left\{e^{{\rm i}kx}:\, k\in \mathbb{Z},\,-N/2\leq k\leq N/2-1 \right\},$$
and
$$S_h=\left\{\phi\in C^{r-2}_p: \phi_{|_{[x_j,x_{j+1}]}} \in\mathbb{P}_{r-1},\, 0\leq j\leq N-1 \right\}.$$
The hybrid spectral-finite element approximation $u_{h}$ of the solution $u$ of (\ref{E21}) is a real $S_{h}$-valued function $u_{h}(t)$ of $t\geq 0$ defined by the ode initial-value problem
\begin{equation}\label{E23}
\begin{split}
&({u_h}_t,\chi)+(\alpha {u_h}_x+\beta u_h {u_h}_x,\chi)+\gamma(P_N\mathcal{G}{u_h}_x,\chi_x)+\delta({u_h}_{xx},\chi_x)=0,\forall \chi\in S_{h}, t\geq 0,\\
&u_h(0)=P_h u_{0},
\end{split}
\end{equation}
where $P_h$, $P_N$ are the $L^2$ projections onto $S_h$ and $S_N$, respectively, given for $w\in L^{2}$ as
$$(P_h w,\chi)=(w,\chi),\quad \forall \chi\in S_h$$
and
$$(P_N w,\phi)=(w,\phi),\quad \forall \phi\in S_N,$$ where $(\cdot,\cdot)$ is the $L^{2}(0,2\pi)$ inner product. For $f\in L^{2}$, $P_{N}f$ is represented by
$$P_N f(x)=\sum_{k=-N/2}^{N/2-1} \hat{f}_k e^{{\rm i}kx},$$
where $\hat{f}_{k}=\frac{1}{2\pi}\int_{0}^{2\pi} f(x)e^{-ikx}dx, k\in \mathbb{Z}$ are the Fourier coefficients of $f$. Note that $\widehat{(\mathcal{G}f)}_{k}=-{\rm i}\sign(k)\hat{f}_{k}$ and that $\mathcal{G}$ is antisymmetric in $L^{2}$.

\subsection{The fully discrete hybrid scheme}
\label{sec22}
We define our fully discrete hybrid scheme following the derivation of the analogous scheme of \cite{BDKMc} in the case of the generalized KdV equation. (This scheme was also used in \cite{BDM}.) Denoting again by $(\cdot,\cdot)$ the $L^{2}(0,2\pi)$ inner product, we define, for each $t\in [0,T]$, the map $F:S_{h}\rightarrow S_{h}$ by the equation
$$(F(u_h),\chi)=-[(\alpha {u_h}_x+\beta u_h {u_h}_x,\chi)+\gamma(P_N\mathcal{G}{u_h}_x,\chi_x)+\delta({u_h}_{xx},\chi_x)], \quad \forall \chi\in S_{h}.$$
Then, the initial-value problem (\ref{E23}) may be written as 
\begin{eqnarray}
{u_h}_t=F(u_h),\quad 0\leq t\leq T,\quad u_{h}(0)=P_{h}u_{0}.\label{E31}
\end{eqnarray}
In addition to $F$ we define the maps $B:S_h\times S_h\rightarrow S_h$, $\Theta_1:S_h\rightarrow S_h$ and $\Theta_2: S_h\rightarrow S_h$ that satisfy for $v,w\in S_h$ and for all $\chi \in S_h$
$$(B(v,w),\chi)=\frac{1}{2}(\beta vw,\chi') =-\frac{1}{2}(\beta(vw)_x,\chi),$$
$$(\Theta_1 v,\chi)=(\alpha v-\delta v_{xx},\chi'),$$
and
$$(\Theta_2 v,\chi)=-(\gamma P_N \mathcal{G} v_x,\chi').$$ If we put
$$F(v,w):=B(v,w)+\Theta_1 v+ \Theta_2 v,$$
we see that
$$F(v):=F(v,v)=B(v)+\Theta_1 v+\Theta_2 v,$$ where $B(v)=B(v,v)$. The initial-value problem (\ref{E31}) is stiff. It is
discretized in the temporal variable by the 2-stage Gauss-Legendre
implicit Runge-Kutta method, which is fourth-order accurate and has good nonlinear stability properties. It corresponds to the Butcher table
\begin{equation*}
\begin{array}{cc|c}
a_{11} & a_{12} & \tau_1\\
a_{21} & a_{22} & \tau_2 \\
\hline
b_1 & b_2 & \\
\end{array}\,\,\,\, =\,\,\,\,
\begin{array}{cc|c} \frac{1}{4} & \frac{1}{4}-\frac{1}{2\sqrt{3}} &
\frac{1}{2}-
\frac{1}{2\sqrt{3}} \\
\frac{1}{4}+ \frac{1}{2\sqrt{3}} & \frac{1}{4} & \frac{1}{2}+
\frac{1}{2\sqrt{3}}\\
\hline
\frac{1}{2} & \frac{1}{2} & \\
\end{array}.
\end{equation*}
The fully discrete scheme is now specified more precisely. Let
$t^n=nk$, $n=0,1,\ldots, M$, where $T=Mk$. We seek $U^n$ approximating $u_{h}(t^{n})$, and  $U^{n,i}$ in $S_h$,
$i=1,2$, as solutions of the system of
nonlinear equations
\begin{equation}\label{E32}
U^{n,i}=U^n+k\sum_{j=1}^2 a_{ij} F(U^{n,j}),\quad i=1,2,\quad 0\leq n\leq M-1,
\end{equation}
and set
\begin{equation}\label{E33}
U^{n+1}=U^n+k\sum_{j=1}^2 b_j F(U^{n,j}), \quad 0\leq n\leq M-1,
\end{equation}
where $U^{0}=u_{h}(0)$. At each time step we solve the nonlinear system (\ref{E32}) using
Newton's method as follows. Given $n \ge 0$, let $U_0^{n,i}\in S_h$, $i= 1, 2$ be an accurate
enough (see below) initial guess for $U^{n,i}$,
the solution of (\ref{E32}). Then the iterates
of Newton's method (called the {\it outer} iterates for reasons that will
become clear presently) $U_j^{n,i}$, $j= 1,2,\ldots$ ($U_j^{n,i}$ approximates $U^{n,i}$) satisfy the
$2\times 2$ block linear system in $S_h\times S_h$,
\begin{eqnarray}
&&\left[\begin{array}{cc}
I+k a_{11} J(U_j^{n,1}) & k a_{12}J(U_j^{n,2})\\
k a_{21}J(U_j^{n,1}) & I+k a_{22} J(U_j^{n,2})
\end{array}\right]\, \left[\begin{array}{c} U_{j+1}^{n,1}\\  U_{j+1}^{n,2} \end{array} \right]=
\left[\begin{array}{c} U^n\\  U^n \end{array} \right]\label{E34}\\
&&-k
\left[\begin{array}{cc}
a_{11} &  a_{12}\\
a_{21} &  a_{22}
\end{array}\right]\, \left[\begin{array}{c} B(U_{j}^{n,1})\\ B( U_{j}^{n,2}) \end{array} \right],\nonumber
\end{eqnarray}

where, for $\psi, \phi$ in $S_{h}$ $$J(\phi)\psi=J_1(\phi)\psi+J_2(\phi)\psi,$$
$$J_1(\phi)\psi=-2B(\phi,\psi)-\Theta_1\psi,$$
and
$$J_2(\phi)\psi=-\Theta_2\psi.$$

The equations (\ref{E34}) represent a $2N\times 2N$ linear system for the coefficients of the new Newton iterates $U_{j+1}^{n,i}$, $i=1,2$, for each $j$, with respect to a basis of $S_{h}$. The two operator equations in (\ref{E34}) are uncoupled as follows: We evaluate the entries of the matrix in the left-hand side of (\ref{E34}) at a point $U^{\ast}\in S_h$, defined by
\begin{equation}\label{E35}
U^{\ast}=\frac{1}{2}(U_0^{n,1}+U_0^{n,2}),
\end{equation}
(which makes the operators in the entries of this matrix independent of $j$ and
allows them to commute with each other). We may then write (\ref{E34}) equivalently as
\begin{equation}\label{E36}
\begin{split}
&\left[\begin{array}{cc}
I+k a_{11} J_1(U^{\ast}) & k a_{12}J_1(U^{\ast})\\
k a_{21}J_1(U^{\ast}) & I+k a_{22} J_1(U^{\ast})
\end{array}\right]\, \left[\begin{array}{c} U_{j+1}^{n,1}\\  U_{j+1}^{n,2} \end{array} \right]=\\
&\left[\begin{array}{c} U^n\\  U^n \end{array} \right]-k
\left[\begin{array}{cc}
a_{11} &  a_{12}\\
a_{21} &  a_{22}
\end{array}\right]\, \left[\begin{array}{c} B(U_{j}^{n,1})\\ B( U_{j}^{n,2}) \end{array} \right]\\
&+k\left[\begin{array}{cc}
a_{11} &  a_{12}\\
a_{21} &  a_{22}
\end{array}\right]\,\left[\begin{array}{cc}
J_1(U^{\ast}) -J(U_j^{n,1}) & 0\\
0 & J_1(U^{\ast}) - J(U_j^{n,2})
\end{array}\right]\,  \left[\begin{array}{c} U_{j+1}^{n,1}\\  U_{j+1}^{n,2} \end{array} \right],
\end{split}\end{equation}
for $j\ge 0$, a form that immediately suggests an iterative scheme for approximating $U_{j+1}^{n,i}$, $i=1,2$. This scheme generates {\it inner} iterates denoted by $U^{n,i,\ell}_{j+1}$ for given $n,i,j$, and $\ell=0,1,2,\ldots$ ($U^{n,i,\ell}_{j+1}$ approximates $U^{n,i}_{j+1}$) that are found recursively from the equations
\begin{equation}\label{E37}
\left[\begin{array}{cc}
I+k a_{11} J_1(U^{\ast}) & k a_{12}J_1(U^{\ast})\\
k a_{21}J_1(U^{\ast}) & I+k a_{22} J_1(U^{\ast})
\end{array}\right]\, \left[\begin{array}{c} U_{j+1}^{n,1,\ell+1}\\  U_{j+1}^{n,2,\ell+1} \end{array} \right]=\\
\left[\begin{array}{c} r_{j+1}^{n,1,\ell}\\  r_{j+1}^{n,2,\ell} \end{array} \right],
\end{equation}
for $\ell\ge 0$, where
$$r_{j+1}^{n,i,\ell}=U^n-k\sum_{m=1}^2a_{im}B(U^{n,m}_j)+k\sum_{m=1}^2a_{im}(J_1(U^{\ast})-J(U^{n,m}_j)) U_{j+1}^{n,m,\ell}.$$
The linear system (\ref{E37}) can be solved efficiently as follows: Since $a_{12}a_{21}<0$, it is possible, upon scaling the matrix on the left-hand side of the system by a diagonal similarity transformation, to write it as
\begin{equation}\label{E38}
\left[\begin{array}{cc}
I+\frac{1}{4} k J_1(U^{\ast}) & k J_1(U^{\ast})/4\sqrt{3}\\
k J_1(U^{\ast})/4\sqrt{3} & I+\frac{1}{4}k  J_1(U^{\ast})
\end{array}\right]\, \left[\begin{array}{c} U_{j+1}^{n,1,\ell+1}\\ \mu  U_{j+1}^{n,2,\ell+1} \end{array} \right]=\\
\left[\begin{array}{c} r_{j+1}^{n,1,\ell}\\  \mu r_{j+1}^{n,2,\ell} \end{array} \right],
\end{equation}
where $\mu=2-\sqrt{3}$. The system (\ref{E38}) is equivalent to the single complex $N\times N$ system
\begin{equation}\label{E39}
(I+k\zeta  J_1(U^{\ast}))Z=R,
\end{equation}
where $\zeta=\frac{1}{4}+{\rm i}/4\sqrt{3}$, and where $Z$ and $R$ are complex-valued functions with real and imaginary parts in $S_h$ which depend upon $n$, $\ell$ and $j$ and are given by
\begin{equation}\label{E310}
Z=U_{j+1}^{n,1,\ell+1}+{\rm i} \mu U_{j+1}^{n,2,\ell+1},\quad R=r_{j+1}^{n,1,\ell+1}+{\rm i}\mu r_{j+1}^{n,2,\ell+1}.
\end{equation}

In practice only a finite number of outer and inner iterates are computed at each time step. Specifically, for $i=1,2$, $n\geq 0$, we compute approximations to the outer iterates $U_j^{n,i}$ for $j=1,\ldots, J_{out},$ for some small positive integer $J_{out}$. For each $j$, $0\leq j\leq J_{out}-1$, $U_{j+1}^{n,i}$ is approximated by the last inner iterate $U_{j+1}^{n,i,J_{inn}}$ of the sequence of inner iterates $U_{j+1}^{n,i,\ell}$, $0\leq \ell \leq J_{inn}$ that satisfy linear systems of the form (\ref{E39}). $J_{inn}$ and $J_{out}$ are such that
$$\left(\sum_{k=1}^2\|U_{j+1}^{n,k,\ell+1}-U_{j+1}^{n,k,\ell}\|^2_{\ell_2}\right)^{1/2}\leq \varepsilon,$$
and
$$\left(\sum_{k=1}^2\|U_{j+1}^{n,k}-U_{j}^{n,k}\|^2_{\ell_2}\right)^{1/2}\leq \varepsilon,$$
where $\|v\|_{\ell_2}$ denotes the Euclidean norm of the coefficients of $v\in S_h$ with respect to its basis, and $\varepsilon$ is usually taken to be $10^{-10}$.

Given $U^n$, the required starting values $U_0^{n,i}$ for the outer (Newton) iteration are computed by extrapolation from previous values as
\begin{equation}\label{E311}
U_0^{n,i}=\alpha_{0,i} U^{n}+\alpha_{1,i} U^{n-1}+\alpha_{2,i} U^{n-2}+\alpha_{3,i} U^{n-3},
\end{equation}
for $i=1,2$, where the coefficients $\alpha_{j,i}$ are such that $U_0^{n,i}$ is the value at $t=t^{n,i}$ of the Lagrange interpolating polynomial of degree at most 3 in $t$ that interpolates to the data $U^{n-j}$ at the four points $t^{n-j}$, $0\leq j\leq 3$. (If $0\leq n\leq 2$, we use the same linear combination, putting $U^j=U^0$ if $j<0$.)

The integrals involving the local terms are computed in general using the 5-point Gauss-Legendre quadrature rule in each spatial interval. The inner product $(P_N\mathcal{G}{u_h}_x,\chi_x)$ involving the nonlocal term is computed  as the inner product $(I_N\mathcal{G}{u_h}_x,\chi_x)$ where the Fourier interpolant $I_N$ is defined as
\begin{eqnarray}
\label{E312}
I_N v(x)=\sum_{k=-N/2}^{N/2-1}\hat{v}_ke^{{\rm i}kx},
\end{eqnarray} where by $\hat{v}_k$ we denote the discrete Fourier coefficients of $v$,
computed  by the Fast Fourier Transform. The inner product $(\cdot,\cdot)$ is approximated by the trapezoidal quadrature rule, which is very accurate for periodic functions.

In the sequel, we shall use the fully discrete scheme described above with the $C^{2}$ cubic splines ($r=4$) as the finite element subspace $S_{h}$. We shall refer to this method as the {\it hybrid scheme/method}.

We checked numerically the orders of convergence of the hybrid scheme as follows. Due to lack of analytical formulas for solutions of the Benjamin equation we considered the nonhomogeneous equation
\begin{equation}\label{E313}
u_t+uu_x+\mathcal{G} u_{xx}+\frac{1}{2}u_{xxx}=f(x,t),\quad (x,t)\in [-1,1]\times [0,T],
\end{equation}
with periodic boundary conditions and $$f(x,t)=e^t\left(\sin(\pi x)+\frac{\pi}{2}e^t\sin(2\pi x)+\left(\pi^2-\frac{\pi^3}{2}\right)\cos(\pi x) \right).$$ The specific equation has a solution $u(x,t)=e^t \sin(\pi x)$.
We solved it numerically up to $T=1$ and we computed the discrete maximum error on the quadrature nodes and the normalized $L^2$ error defined as $\|e_h(\cdot,t^n)\|/\|e_h(\cdot,0)\|,$ where $e_{h}=u-U$. The numerical method appears to converge with an optimal rate in space ($r=4$) but with a suboptimal rate equal to three in time.
\begin{table}[ht]
\centering
\begin{tabular}{cccccc}
\hline
 $N$ & $M$ & $L^{\infty}$ Error & Rate & $L^2$ Error & Rate\\
 \hline
   $4$ & $1000$ & $ 0.2630\times 10^{-1}$ & -- & $0.4263\times 10^{-1}$ & --\\
   $8$ & $1000$ &  $0.2654\times 10^{-2}$ & $3.309$ & $0.4125\times 10^{-2}$ & $3.370$\\
  $16$ & $1000$ & $ 0.1916\times 10^{-3}$ & $3.793$ & $0.2686\times 10^{-3}$ & $3.941$\\
  $32$ & $1000$ &  $0.1243\times 10^{-4}$ & $3.945$ & $0.1693\times 10^{-4}$ & $3.988$\\
  $64$ & $1000$ &  $0.7863\times 10^{-6}$ & $3.983$ & $0.1060\times 10^{-5}$ & $3.997$\\
 $128$ & $1000$ & $ 0.5068\times 10^{-7}$ & $3.956$ & $0.6636\times 10^{-7}$ & $3.998$\\
 \hline
\end{tabular}
\caption{\label{tav1}Spatial rates of convergence (hybrid scheme)}
\end{table}

\begin{table}[ht]
\centering
\begin{tabular}{cccccc}
\hline
 $N$ & $M$ & $L^{\infty}$ Error & Rate & $L^2$ Error & Rate\\
 \hline
  $20$ & $ 20$ &  $0.1301\times 10^{-3}$ & -- & $0.1249\times 10^{-3}$ & --\\
  $40$ &  $40$ &  $0.1866\times 10^{-4}$ & $2.802$ & $0.1678\times 10^{-4}$ & $2.896$\\
  $80$ &  $80$ &  $0.3888\times 10^{-5}$ & $2.262$ & $0.3733\times 10^{-5}$ & $2.169$\\
 $160$ & $160$ &  $0.5566\times 10^{-6}$ & $2.804$ & $0.5465\times 10^{-6}$ & $2.772$\\
 $320$ & $320$ &  $0.7289\times 10^{-7}$ & $2.933$ & $0.7101\times 10^{-7}$ & $2.944$\\
 $640$ & $640$ &  $0.9443\times 10^{-8}$ & $2.948$ & $0.8994\times 10^{-8 }$& $2.981$\\
 \hline
\end{tabular}
\caption{\label{tav2}Temporal rates of convergence (hybrid scheme)}
\end{table}
Tables \ref{tav1} and \ref{tav2} show the numerical spatial and temporal rates of convergence of the error for this experiment computed in the discrete maximum norm and the normalized $L^{2}$ norm at $t=T=1$. Here $N$ is the number of spatial intervals and $M=T/k$. We observe that the spatial rate is practically optimal (four) and that the temporal rate approximates the value $p=3$ as $N, M$ increase. (For this experiment, with the tolerance set at $\epsilon=10^{-10}$, the number of Newton iterations $J_{out}$ came out to be always one and $J_{inn}$ varied in general between one and four provided $k$ and $h$ were sufficiently small.) The theoretical order of accuracy of the two-stage Gauss-Legendre RK method is of course equal to four and this value is observed experimentally for the KdV equation, i.~e. when the nonlocal term $\mathcal{G}u_{xx}$ is not present, see e.~g. (\cite{BDKMc}, Table 3). In our case, the loss of one order of temporal accuracy is apparently caused by the presence of the nonlocal term: Observe that in the Jacobian $J_{1}(U^{*})$ in the matrix of operators in the left-hand side of (\ref{E36}) we did not include the part of the Jacobian $J_{2}=-\Theta_{2}$ corresponding to the nonlocal term but transferred it to the right-hand side, in order to retain sparsity in the operators on the left when a basis of small support is chosen for $S_{h}$. This efficiency consideration renders the scheme explicit with respect to the nonlocal term and linearly implicit with respect to the rest of the terms in the equation, and causes the loss of temporal accuracy by one order.

We did not detect any need for a stability bound on $k/h$ for these computations. (Values as high as $k/h=8$ were tried.) Of course accuracy is reduced as $k$ increases and so in the numerical experiments of sections \ref{sec3}-\ref{sec5} $k/h$ was taken much smaller.

In the sequel, we shall also on occasion compute solutions of the Benjamin-Ono (BO) equation, mainly in order to test our numerical schemes. (BO is a good testing ground for our purposes since it has solitary-wave solutions that are known in closed form and are not trivial to simulate on a finite interval as they decay like $O(x^{-2})$ as $|x|\rightarrow \infty$. In addition, their interactions are \lq clean\rq\ due to the integrability of the BO.) For this reason, we briefly report on the performance of the hybrid method in the case of the BO. It is easy to verify, to begin with, that the spatial rate of convergence is again equal to $4$. However, we found that the explicit way that the Newton solver treats the nonlocal term causes the hybrid method to converge under a stability condition of the form $k=\alpha h^{2}$. (In the case of the example (\ref{E313}) with no KdV term, $\alpha \cong 0.6$ was sufficient.)

In the case of the Benjamin-Ono equation, due to the restrictive stability condition $k=\alpha h^2$, if we take a fixed number $N$ of spatial intervals, we observe that the errors cease to decrease at a certain point because the temporal error becomes much smaller than the spatial error. It is thus not easy to compute the asymptotic rate of the temporal error. To accomplish this we did the following: For a fixed value of $h$, we solved the problem in the domain $[-15,15]$ with the hybrid method up to $T=1$ for various values of $k$. We chose $h=0.05$ (i.e. $N=600$) to ensure that the spatial errors will be larger than the temporal errors. We also chose a reference value of $k=k_{ref}=10^{-4}$ ($M=10000$) and we computed the solution $U_{ref}$. We then chose values of $k$ larger than $k_{ref}$ but small enough so as to satisfy the stability condition and computed $U_k$ and the normalized errors
$$E^{\ast}(T)=\frac{\|U_{ref}(T)-U_k(T)\|}{\|u(0)\|}.$$

It turns out that for small values of $k$, which are nevertheless considerably larger than $k_{ref}$, the expected temporal rate of convergence is visible because subtracting $U_{ref}(T)$ from $U_k(T)$, essentially cancels the spatial error of the latter approximation. The results of these computations are presented in Table \ref{tav3}.

\begin{table}[ht]
\centering
\begin{tabular}{cccccc}
\hline
 $N$ & $M$ & $L^{\infty}$ Error & Rate & $L^2$ Error & Rate\\
 \hline
 $600$ & $1250$ & $0.7454\times 10^{-7}$ &--&  $0.7947\times 10^{-7}$ &--\\
 $600$ & $1600$ & $0.3528\times 10^{-7}$  & $3.030$ & $0.3783\times 10^{-7}$ & $3.007$\\
 $600$ & $2000$  & $0.1797\times 10^{-7}$ & $3.024$ & $0.1930\times 10^{-7}$ & $3.016$\\
 $600$ & $2500$ &  $0.9165\times 10^{-8}$  & $3.018$ & $0.9808\times 10^{-8}$  & $3.033$\\
 $600$ & $3200$ &  $0.4298\times 10^{-8}$ & $3.068$ & $0.4585\times 10^{-8}$ & $3.081$\\
 $600$ & $4000$ & $0.2129\times 10^{-8}$  & $3.148$ & $0.2272\times 10^{-8}$ & $3.146$\\
 \hline
\end{tabular}
\caption{\label{tav3}Temporal rates of convergence for BO (hybrid scheme)}
\end{table}
\subsection{A fully discrete pseudospectral scheme}
\label{sec23}
In addition to the hybrid method, we shall use for checking purposes a spectral method. For continuous $2\pi-$periodic complex-valued functions $u, v$ we let $(u,v)_N:=\frac{2\pi}{N}\sum_{j=0}^{N-1}u(x_j)\overline{v(x_j)}$. We consider the following semidiscrete Fourier-collocation (pseudospectral) scheme, cf. \cite{M,CHQZ}, that approximates the solution $u$ of (\ref{E21}) on $[0,2\pi]$ by $u^{N}\in S_{N}$ defined by the equations
\begin{equation}\label{E314}
\begin{array}{l}
(u^N_t+[\alpha u^N+(\beta/2) (u^N)^2-\gamma {\mathcal G} u^N-\delta u^N_{xx}]_x,\chi)_N=0,\quad \forall \chi\in S_{N}, t\geq 0,\\
u^N(x,0)=I_N u_0,
\end{array}
\end{equation}
where $I_{N}$ is given by (\ref{E312}).
By choosing $\chi=e^{-{\rm i}kx}$ for $k=-N/2,\ldots,N/2-1$, we obtain the following system of ode's for the Fourier coefficients $\hat{u}_{k}$ of $u^N$ for $k=-N/2,\ldots,N/2-1$:

\begin{equation}\label{E315}
\frac{d}{dt}\hat{u}_{k}+\frac{\beta}{2}{\rm i}k(\hat{u}\ast\hat{u})_k+\omega(k)\hat{u}_k=0,\quad t\geq 0,\quad
\hat{u}_{k}(0)=\widehat{I_{N}u_{0}}_{k},
\end{equation}
where
$$\omega(k)=\alpha {\rm i}k-\gamma {\rm i} |k| k+\delta {\rm i} k^3.$$
Multiplying  the ode's by $e^{\omega(k)t}$ and
setting $\hat{U}_k=e^{\omega(k)t} \hat{u}_k$ we may write them as
\begin{equation}\label{E316}
\frac{d}{dt}\hat{U}_{k}+\frac{\beta}{2}{\rm i}ke^{\omega(k)t}\left[(e^{-\omega(k)t}\hat{U})\ast(e^{-\omega(k)t}\hat{U})\right]_k=0.
\end{equation}
To compute the convolution $\ast$ we use the formula ${\mathcal F}([{\mathcal F^{-1}}(e^{-\omega(k)t}\hat{U})]^2)$,
where ${\mathcal F}$ is the discrete Fourier transform.
The resulting ode system is discretized by the explicit classical fourth-order Runge-Kutta method in time. Hence, this fully discrete scheme belongs to the class of the so-called \lq integrating factor\rq\ schemes, \cite{CK,MT,KT}, having improved stability properties, as they attempt to reduce stiffness. (The last-quoted paper has a useful review of related schemes.)

We verified the fourth order of temporal accuracy of this scheme by computing its errors in the case of the nonhomogeneous problem (\ref{E313}) at $t=1$ for $N=100$ and an increasing number of time steps. The results are shown in Table \ref{tav4}. (The numerical temporal rate in the case of the analogous numerical experiments for the BO equation was also found to be $4$.)

\begin{table}[ht]
\centering
\begin{tabular}{cccccc}
\hline
 $N$ & $M$ & $L^{\infty}$ Error & Rate & $L^2$ Error & Rate\\
 \hline
 $100$ &  $400$ &  $0.1695\times 10^{-7}$ & --      &     $0.6240\times 10^{-8}$ &  --\\
 $100$ &  $800$ &  $0.1082\times 10^{-8}$ & $3.969$ & $0.3900\times 10^{-9}$ & $4.000$\\
 $100$ & $1600$  & $0.6839\times 10^{-10}$ & $3.984$ & $0.2437\times 10^{-10}$ & $4.000$\\
 $100$ & $3200$  & $0.4305\times 10^{-11}$ & $3.990$ & $0.1526\times 10^{-11}$ & $3.998$\\
 $100$ & $6400$  & $0.2718\times 10^{-12}$ & $3.986$ & $0.9494\times 10^{-13}$ & $4.006$\\
 \hline
\end{tabular}
\caption{\label{tav4}Temporal rates of convergence (spectral scheme).}
\end{table}
We shall henceforth refer to this fully discrete pseudospectral scheme as the \lq spectral\rq\ method.

%NUMERICAL EXPERIMENTS FOR BO: FOR EXAMPLE THE TABLE IN THE BO CASE FOR THE SAME EXACT SOLUTION BUT A SLIGHTLY DIFFERENT RHS IS THE FOLLOWING TABLE \ref{tav5}. (NOTE: DIMITRI SUGGESTS NOT TO ADD IT).
%
%\begin{table}[ht]\caption{\label{tav5}Time rates of convergence (Fourier collocation) for BO}
%\begin{tabular}{cccccc}
% N & M & MAXERR & RATE & L2ERR & RATE\\
% \hline
% 100 & 400  & 0.1705D-07 & 0.000 & 0.5540D-08 & 0.000\\
% 100 & 800  & 0.1066D-08 & 3.999 & 0.3462D-09 & 4.000\\
% 100& 1600  & 0.6665D-10 & 4.000 & 0.2164D-10 & 4.000\\
% 100 &3200  & 0.4168D-11 & 3.999 & 0.1356D-11 & 3.997\\
% 100 &6400  & 0.2465D-12 & 4.080 & 0.7910D-13 & 4.099\\
% \hline
%\end{tabular}
%\end{table}

\subsection{Validation of the hybrid method}
\label{24}
We now present the results of some numerical tests that we performed with both schemes in order to validate further the hybrid method and compare its results with those of the spectral scheme.

\begin{figure}[!htbp]
\centering
{\includegraphics[width=\textwidth]{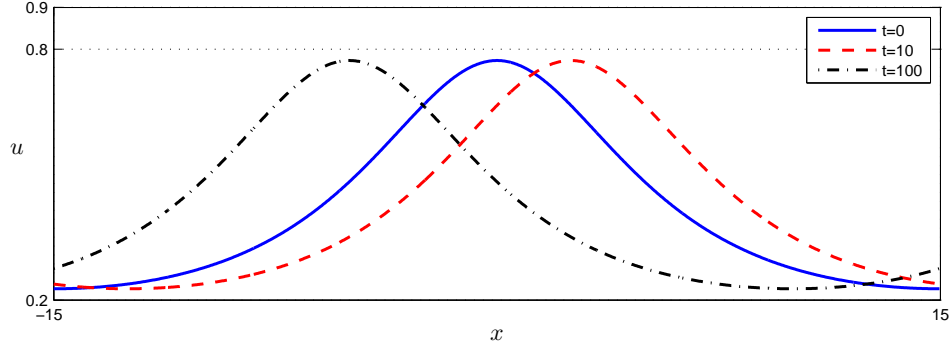}}
\caption{Numerical evolution of the periodic-travelling wave solution (\ref{E317}) of the Benjamin-Ono equation.}\label{F1}
\end{figure}

In our first experiment we simulate the propagation of a periodic travelling-wave solution of the Benjamin-Ono equation that was used in \cite{TM}. This solution resembles a solitary wave and is given by the formula
\begin{equation}\label{E317}
u(x,t)=\frac{2c_s A^2}{1-\sqrt{1-A^2}\cos (c_s A(x-c_st))},
\end{equation}
where $A=\frac{\pi}{c_s L}$. This is a $2L-$periodic solution of the BO with coefficients $\alpha=\delta=0$, $\beta=\gamma=1$ in (\ref{E11}). We approximated it by the spectral method with $N=1024, k=0.02$ and the hybrid method in two runs with $N=256$ and $k=0.01$ and with $N=1024$ and $k=5\times 10^{-4}$, respectively, on the interval $[-L,L]$ with $L=15$ and $c_{s}=0.25$ for $0\leq t\leq 100$, using (\ref{E317}) at $t=0$ as initial condition. The numerical solution is shown in Figure \ref{F1} at $t=0, 10$ and $100$. (All three numerical profiles coincided within graph thickness.)

\begin{figure}[!htbp]
\centering
{\includegraphics[width=\textwidth]{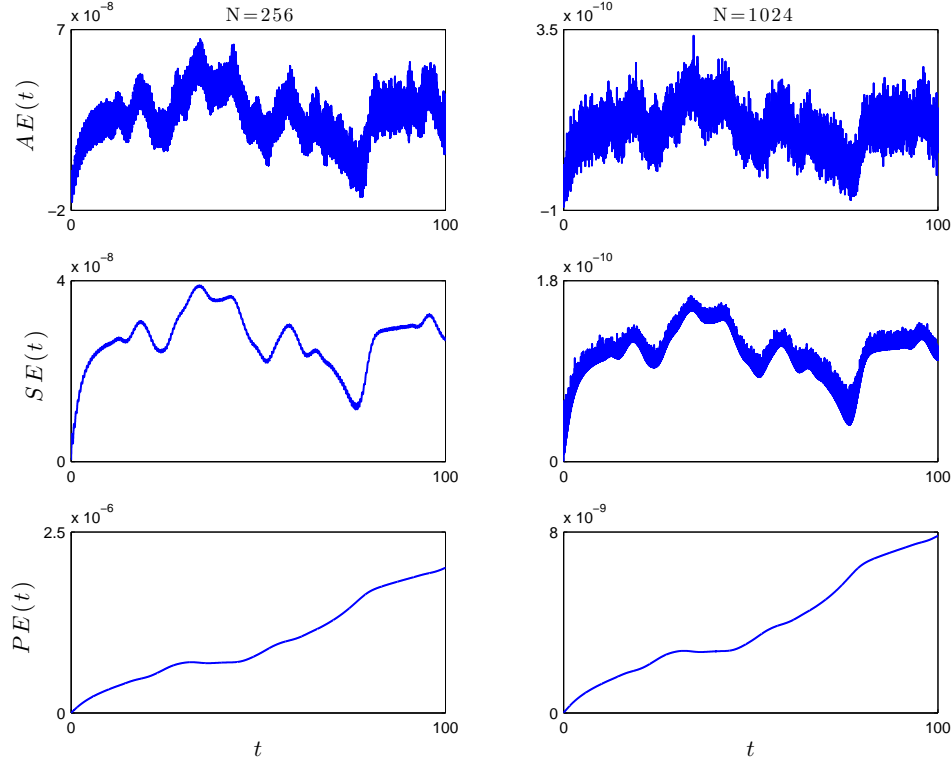}}
\caption{Amplitude ($AE(t^{n})$), Shape ($SE(t^{n})$) and Phase ($PE(t^{n})$) errors of the hybrid scheme for $N=256, 1024$, approximating the solution (\ref{E317}) of the BO equation}
\label{F2}
\end{figure}

In this example, the errors of the spectral method were all in the range $10^{-9}$ to $10^{-11}$. In the two runs of the hybrid scheme, the normalized $L^{2}$ error, defined as $\max_{n}\frac{||u(t^{n})-U^{n}||}{||U^{0}||}$, was of $O(10^{-7})$ for $N=256$ and of $O(10^{-11})$ for $N=1024$. In both cases, the $L^{2}$ norm of the numerical solution was equal to $2.50662827463$ while the Hamiltonian (invariant $E(u)$ given by (\ref{E14})) was equal to $-0.473444593881$. (Both were preserved for $0\leq t\leq 100$ up to the twelve significant digits shown.) In addition, for the hybrid scheme we computed for each $t^{n}$ several other types of errors that are relevant in assessing the accuracy of approximation of solitary-type waves, cf. \cite{BDKMc, BDM}. These were: (i) The {\it (normalized) amplitude error} $AE(t^{n})=\left|\frac{u_{max}-U^{n}(x^{*})}{u_{max}}\right|$, where $u_{max}$ is the maximum value of the exact solution and $x^{*}$ is the point where the approximate solution $U^{n}$ achieves its maximum, found by applying Newton's method to compute the root of the equation $\frac{d}{dx}U^{n}(x)=0$ that corresponds to the maximum of $U^{n}$. (ii) The $L^{2}$ {\it (normalized) shape error} defined as $SE(t^{n})=\inf_{\tau}||U^{n}-u(\cdot,\tau)||/||u_{0}||$, computed as $SE(t^{n})=\xi(\tau^{*})$, where $\tau^{*}$ is the point near $t^{n}$ (found by Newton's method) where $\frac{d}{d\tau}(\xi^{2})=0$, with $\xi(\tau)=||U^{n}-u(\cdot,\tau)||/||u_{0}||$. (iii) The associated {\it phase error} $PE(t^{n})=\tau^{*}-t^{n}$. Figure \ref{F2} shows these errors as functions of $t^{n}$ up to $T=100$, for $N=256$ and $N=1024$. The speed $c_{s}=0.25$ of the travelling wave was preserved for $N=256$ to $6$ digits up to $t=50$ and to $5$ digits up to $t=100$, while for $N=1024$ up to at least $7$ digits up to $t=100$.

In a second validation experiment we computed the evolution of a solitary wave for the Benjamin equation (\ref{E21}) with $\gamma=0.5$ (all other coefficients being equal to one) with $L=128$ up to $T=100$. The initial solitary-wave profile was generated with high accuracy by numerical continuation with the CGN method as explained in \cite{DDM1} and in Section \ref{sec3} of the present paper. We solved the problem by the hybrid and the spectral schemes. Table \ref{tav5b} presents the results of two runs with comparable errors for this problem. The spectral method is faster by a factor
\begin{table}[ht]
\centering
\begin{tabular}{ccc}
\hline
  & Hybrid & Spectral \\
 \hline
 $N$ & $2048$  & $256$\\
 $k$ & $1\times 10^{-2}$  & $1\times 10^{-2}$\\
 $L^{2}$ error& $0.4398\times 10^{-6}$  & $0.8024\times 10^{-6}$\\
 $H^{1}$ error & $0.3664\times 10^{-6}$ & $0.8888\times 10^{-6}$\\
 $SE$ & $0.1370\times 10^{-6}$ & $0.1117\times 10^{-5}$\\
 $PE$ & $0.1728\times 10^{-5}$ &$0.4642\times 10^{-7}$\\
 $H$ & $0.4827201809$ & $0.482720$\\
 cpu time (sec)& $59$ & $30$ \\
 \hline
\end{tabular}
\caption{\label{tav5b} Errors at $T=100$ and parameters for the hybrid and spectral methods. Solitary wave, Benjamin equation, $\gamma=0.5$}
\end{table}
 of two but the hybrid method conserves the Hamiltonian $H=I+E$ up to $10$ digits, four more than in the case of the spectral method. In the table the $L^{2}$ and shape errors are normalized as explained earlier. The (normalized) $H^{1}$ error, defined analogously, is a useful error metric for oscillatory profiles such as the solitary waves of the Benjamin equation.

\begin{figure}[!htbp]
\centering
{\includegraphics[width=\textwidth]{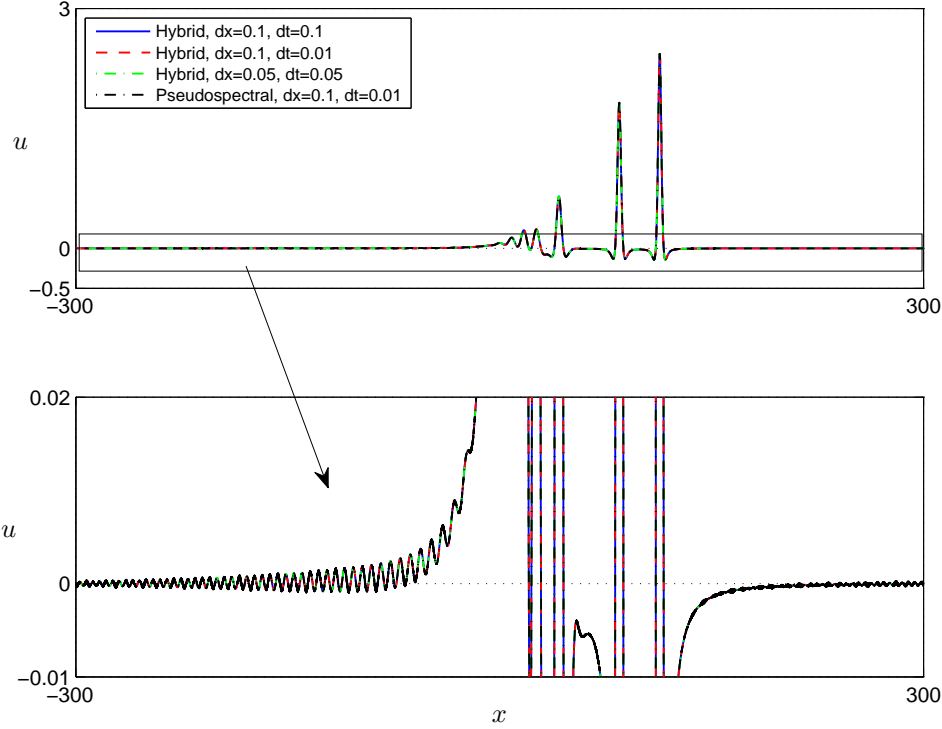}}
\caption{Resolution of the \lq Gaussian\rq\ $2e^{-(x/4)^{2}}$ into solitary waves. Benjamin equation, T=100. The profile on the bottom is a magnification of that on the top.}\label{F3}
\end{figure}

 In our third experiment we solved the Benjamin equation in the form $u_{t}+uu_{x}+\mathcal{G}u_{xx}+u_{xxx}=0$ for $x\in [-300,300]$ up to $T=100$ using as initial condition the Gaussian $u(x,0)=2e^{-(x/4)^{2}}$. As expected, \cite{KB}, the initial profile resolves itself into a series of solitary waves. As Figure \ref{F3} shows, by $T=100$ three solitary waves have appeared, followed by a dispersive tail.

We used the solution obtained by the spectral scheme with $N=6000, k=0.01$ as the benchmark and recomputed the solution with the hybrid scheme for various values of the discretization parameters $h$ and $k$ starting from $h=0.1, k=0.1$ and reducing $h$ and/or $k$. Some of the profiles produced by the hybrid runs are shown in Figure \ref{F3}; they all coincide within graph thickness with the spectral solution. (It should be mentioned that the spectral scheme with $k=600/N$ blew up and needed $k=O((600/N)^{2})$ for stability.)
\section{Generation and propagation of solitary waves}
\label{sec3}
In this section we first review the numerical technique that we used to generate solitary-wave solutions of the Benjamin equation. These solitary-wave profiles were taken as initial values for the hybrid time-stepping method and integrated forward in time. We present in some detail the temporal evolution of various error metrics suitable for assessing the accuracy of these numerically generated travelling waves.

As was already mentioned in the Introduction, the solitary waves of the Benjamin equation  are travelling-wave solutions of (\ref{E11}) of the form $u(x,t)=\varphi(x-c_{s}t), c_{s}>0$, such that $\varphi$ and its derivatives tend to zero as $\xi=x-c_{s}t$ approaches $\pm \infty$. Consequently, $\varphi$ satisfies the equation (\ref{E15}), from which, taking Fourier transforms, we obtain
\begin{eqnarray*}
(-c_{s}+\alpha-\gamma |k|+\delta k^{2})\widehat{\varphi}
+\frac{\beta}{2}\widehat{\varphi^{2}}=0,\quad k\in\mathbb{R},
\end{eqnarray*}
where $\widehat{\varphi}(k)$ is the Fourier transform of $\varphi$. If we discretize this equation assuming periodic boundary conditions on $[-L,L]$ and using the discrete Fourier transform to compute the convolution as in section \ref{sec23}, we obtain the $N\times N$ nonlinear system of equations
\begin{eqnarray}\label{E41}
(-c_{s}+\alpha-\gamma |k|+\delta k^{2})\widehat{\varphi^{N}}_{k}
+\frac{\beta}{2}\left(\widehat{\varphi^{N}\ast \varphi^{N}}\right)_{k}=0,\quad k=-\frac{N}{2},\ldots,\frac{N}{2}-1,
\end{eqnarray}
where $\varphi^{N}$ is the approximation of $\varphi$ in $S_{N}$ and $\widehat{\varphi^{N}}_{k}$ denotes its $k^{\rm th}$ Fourier coefficient.

To solve (\ref{E41}) we use an incremental continuation technique with respect to the parameter $\gamma$, following e.~g. \cite{ABR}. For a fixed set of constants $\alpha, \beta, \delta, c_{s}$ in (\ref{E41}) we consider a homotopic path $\gamma_{0}=0<\gamma_{1}<\ldots<\gamma_{M}=\gamma$ and solve (\ref{E41}) successively for $\gamma_{0},\gamma_{1},\ldots,\gamma_{M}$ with an iterative nonlinear solver, using for each $j$ the numerical solution for $\gamma=\gamma_{j-1}$ as an initial guess in solving for $\gamma=\gamma_{j}$. (The starting value $\gamma_{0}=0$ of the path corresponds to the KdV equation for which exact solitary-wave solutions are available.) The incremental continuation technique has the added advantage that it produces a series of solitary waves for varying values of $\gamma$ with a fixed speed $c_{s}$.

The nonlinear system solver that we used to generate the solution of (\ref{E41}) for each $\gamma_{j}$ was Newton's method, wherein the attendant linear systems were solved by an inner iteration performed by the preconditioned conjugate gradient technique. The resulting iterative scheme, called CGN in the sequel, was described in detail in \cite{DDM1}, where it was also compared with several other nonlinear solvers and found to be more efficient, with respect to a variety of metrics, for approximating solutions of (\ref{E41}). We refer the reader to \cite{DDM1} for the implementation of CGN; let us just mention that for the computations in the present paper the Newton iteration was terminated when the quantity $||\varphi_{[\nu]}^{N}-\varphi_{[\nu -1]}^{N}||/||\varphi_{[\nu]}^{N}||$ became less than $10^{-15}$. (Here $\varphi_{[\nu]}^{N}$ is the $\nu$-th Newton iterate approximating $\varphi^{N}$). The preconditioned conjugate-gradient inner iteration was terminated when $||R^{(i)}||_{M}/||R^{(0)}||_{M}$ became less than $10^{-2}$. Here $R^{(i)}$ is the residual defined in the standard way in the conjugate-gradient algorithm, and the norm $||\cdot||_{M}$ is the weighted $L^{2}$ norm $(\cdot,M^{-1}\cdot)^{1/2}$, where $M=cI-\partial_{xx}$ is the preconditioning operator that we used; its action in Fourier variables is $c+k^{2}$ and the value $c=0.275$ was found to be optimal in computations. The number of CG inner iterations needed to reach the threshold defined above varied between $3$ and $10$ typically.
 \begin{figure}[!htbp]
\centering
{\includegraphics[width=\textwidth]{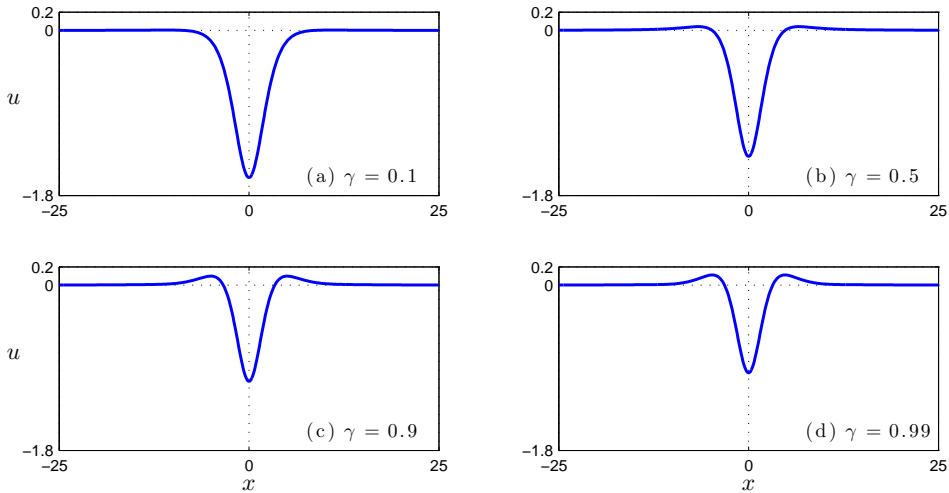}}
\caption{Solitary waves of the Benjamin equation for various values of $\gamma$, $c_s=0.45$.}%
\label{F4}
\end{figure}
 \begin{figure}[!htbp]
\centering
{\includegraphics[width=\textwidth]{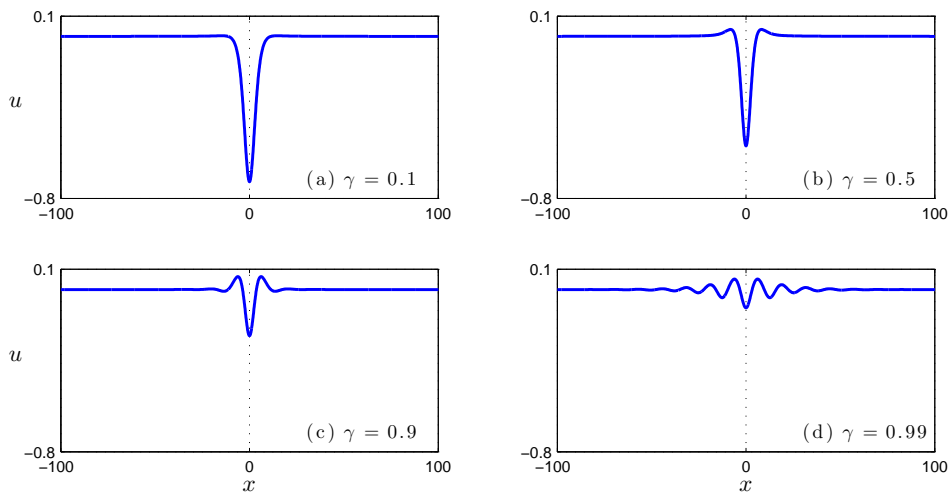}}
\caption{Solitary waves of the Benjamin equation for various values of $\gamma$, $c_s=0.75$.}%
\label{F5}
\end{figure}

Using this algorithm we produced solitary waves of the Benjamin equation in $[-256,256]$ with $N=4096$ using $\gamma_{j}=j\Delta\gamma, j=1,\ldots,99$, with $\Delta \gamma=0.01$ and an exact solitary wave of the KdV equation at $\gamma_{0}=0$. In all computations we took $\alpha=\beta=\delta=1$. Figure \ref{F4} shows the computed profiles of the solitary waves for $c_{s}=0.45$ and $\gamma=0,1,0.5,0.9,0.99$, while Figure \ref{F5} shows the solitary waves corresponding to $c_{s}=0.75$ for the same values of $\gamma$. As is well-known, the number of oscillations increases with $c_{s}$ and $\gamma$.

\begin{figure}[!htbp]
\centering
{\includegraphics[width=\textwidth]{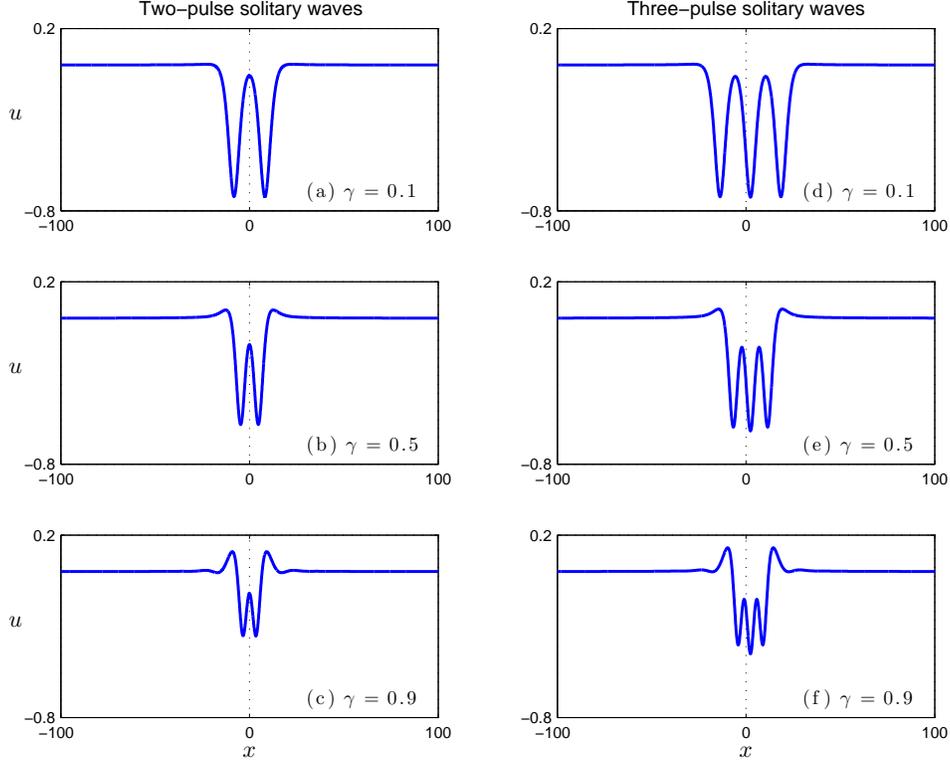}}
\caption{Two-pulse (a,b,c) and three-pulse (b,d,f) solitary waves of the Benjamin equation for $\gamma=0.1,0.5,0.9, c_s=0.75$}%
\label{F6}
\end{figure}

We also constructed with the same technique {\it multi-pulse} solitary waves by starting at $\gamma_{0}=0$ with a superposition of translated KdV solitary waves as explained in \cite{DDM1}. Two-- and three--pulse such solitary waves are shown for $\gamma=0.1,0.5,$ and $0.9$ and $c_{s}=0.75$ in Figure \ref{F6}.

As a measure of the accuracy of the CGN method for approximating the solution of (\ref{E41}) for each value of $\gamma$ we computed the $L^{2}$ norm of the residual $r$, whose $k$-th Fourier component is defined as the left-hand side of (\ref{E41}) with $\phi^{N}$ replaced by its numerical approximation. The value of $||r||$ for single-- and two-- and three-- pulse solitary waves as a function of $\gamma$  remained smaller than $5\times 10^{-13}$ but in general the residual increases as $\gamma$ approaches one, a fact that reflects the difficulty in solving the nonlinear systems with $\gamma$ close to one.

The above-described technique for generating solitary waves of the Benjamin equation  was found to be more accurate, compared to iterative \lq cleaning\rq\ , cf. e.~g. \cite{KB}, wherein one isolates  and \lq cleans\rq\ iteratively solitary waves that are produced by resolution of suitable initial data, and which works well in case the solitary waves decay exponentially. In the case of the Benjamin equation, for which the solitary waves are known to decay quadratically, \cite{B2,CB}, we found that even for large spatial computational intervals it was very hard to make the values at the boundaries of the solitary waves produced by iterative cleaning less than $O(10^{-5})$. This small truncation error produced dispersive oscillations of the same order of magnitude that very fast polluted the ensuing solution when such solitary-wave profiles were used as initial values in evolution studies. Of course, for solitary waves produced by iterative cleaning one does not have {\it a priori} knowledge of their speed, so it is not easy to design systematic experiments with families of solitary waves of varying speed.

We used the numerical solitary waves that we constructed as initial values $u_{0}$ and integrated in time the Benjamin equation using the fully discrete hybrid scheme implemented as in Section \ref{sec2}. As a further test of the accuracy of the numerical solitary waves and the time-stepping technique we computed several invariants of the evolution and various pertinent error measures. In all cases we used the spatial interval $[-256,256]$ and $N=4096$ and we integrated the equation up to $T=300$.

Table \ref{tav5c} shows the values of the $L^{2}$ norm, of the invariant $H=I+E$, where $I$ and $E$ are discrete versions of the quantities defined in (\ref{E13}) and (\ref{E14}), respectively, and of the amplitude of the numerically propagated single-pulse solitary waves with $c_{s}=0.75$ for various values of $\gamma$. The digits shown for each quantity were conserved up to $T=300$.
\begin{table}[ht]
\begin{center}
\begin{tabular}{ccccc}
\hline
$\gamma$ & $L^2$-norm & $H$ & amplitude  \\
\hline
$0.1$ & $1.6096361661$ & $1.09624383030$ &  $-0.7183404$ \\
$0.5$ & $1.08290587306$ &  $0.48258984490$ & $-0.541174$ \\
$0.9$ &  $0.44162186544$ &  $0.07565402212$ & $-0.2280941$ \\
$0.95$ & $0.33588124247$ & $0.04319622837$ & $-0.165667$ \\
$0.99$ & $0.2429264136$ & $0.022247817281$ & $-0.090357$ \\
\hline
\end{tabular}
\end{center}
\caption{Conserved quantities for numerical evolution up to $T=300$ of single-pulse solitary waves of speed $c_s=0.75$ for various values of $\gamma$.}\label{tav5c}
\end{table}

Table \ref{tav6} shows the conserved digits of the same quantities for the analogous propagation experiment with two- and three-pulse solitary waves with $\gamma=0.5$.

\begin{table}[ht]
\begin{center}
\begin{tabular}{cccccc}
\hline
Number of pulses &$\gamma$ & $L^2$-norm & $H$ & amplitude\\
\hline
$2$ & $0.5$ & $1.6419433913$ & $1.1164182800$ & $-0.582995$ \\
$3$ & $0.5$ & $2.0816580537$ & $1.800497679$ & $-0.618111$ \\
\hline
\end{tabular}
\end{center}
\caption{Conserved quantities for numerical evolution up to $T=300$ of multi-pulse solitary waves of speed $c_s=0.75$ for $\gamma=0.5$}\label{tav6}
\end{table}

In these computations the quantity $H$ was defined at $t^{n}$ as
$$\frac{1}{2}\int_{-L}^{L}\left(U^2+\frac{\beta}{3}U^3+\delta U_x^2-\gamma U I_N \mathcal{G} U_x\right)dx,$$ where $U=U^{n}$, the integrals being evaluated by numerical quadrature as described in Section \ref{sec3}.

\begin{figure}[!htbp]
\centering
{\includegraphics[width=\textwidth]{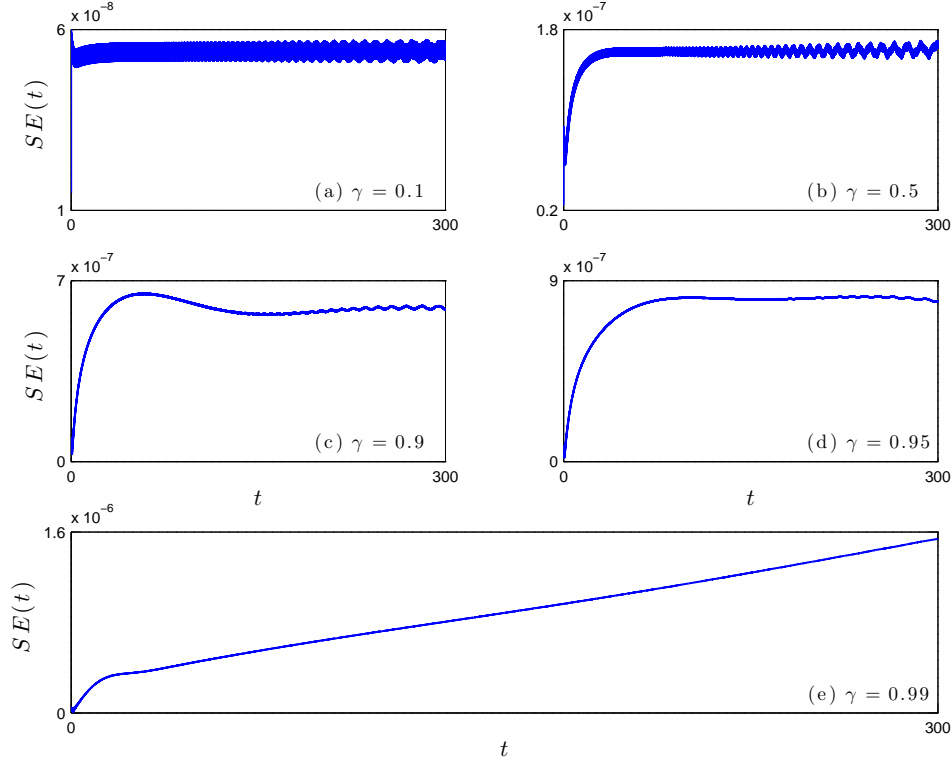}}
\caption{Shape error of the numerical propagation of single-pulse solitary waves with $c_s=0.75$ and various values of $\gamma$.}%
\label{F8}
\end{figure}

\begin{figure}[!htbp]
\centering
{\includegraphics[width=\textwidth]{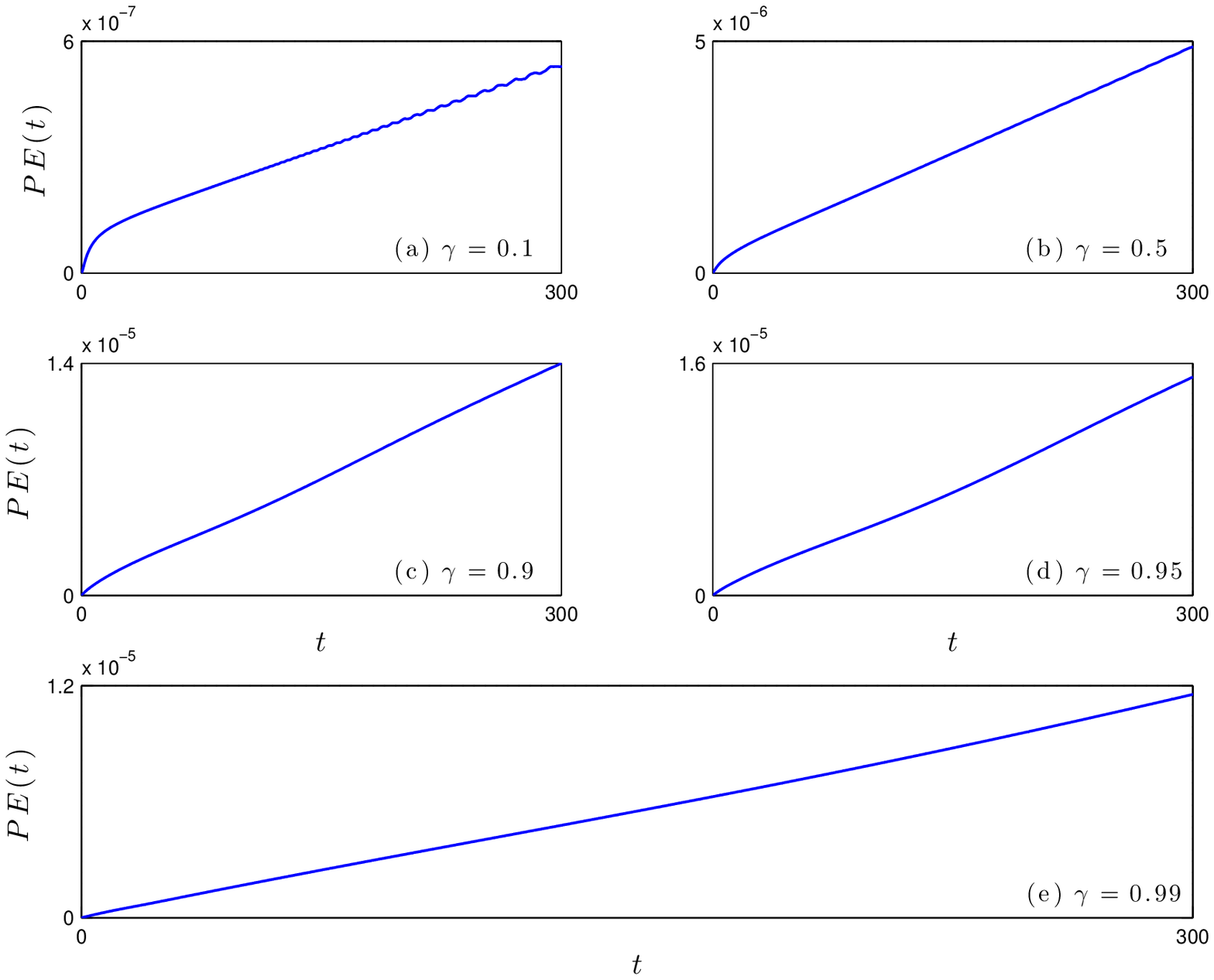}}
\caption{Phase error of the numerical propagation of single-pulse solitary waves with $c_s=0.75$ and various values of $\gamma$.}%
\label{F9}
\end{figure}

In Figure \ref{F8} we show the $L^{2}$ (normalized) shape error of the propagating numerical single-pulse solitary wave for $c_{s}=0.75$ and various values of $\gamma$, as function of $t^{n}$. This quantity is defined as $$SE(t^n)=\inf_{\tau}\|U^n-\varphi_h(\cdot-c_{s}\tau)\|/\|\varphi_h\|,$$ where $\varphi_{h}=P_{h}\varphi^{N}=U^{0}$ is the $L^{2}$-projection on $S_{h}$ of the numerically generated initial solitary wave $\varphi^{N}$.
As in section \ref{sec2}, $SE(t^n)$ is again computed as $\xi(\tau^{*})$, where $\tau^{\ast}$ is the point near $t^{n}$ (found by Newton's method) where $\frac{d}{d\tau}\xi^2(\tau^\ast)=0$, with
$\xi(\tau):=\|U^n-\varphi_h(\cdot-c_{s}\tau)\|/\|\varphi_h\|.$ The shape errors increase with $\gamma$ and stabilize with $t$ except in the case $\gamma=0.99$ where a linear temporal growth is observed. (They range from $O(10^{-8})$ to $O(10^{-6})$.) Figure \ref{F9} shows the analogous graphs for the phase error, defined as $PE(t^n)=\tau^\ast-t^n$. The phase errors increase linearly with $t$ and with $\gamma$ for fixed $t$, ranging from $O(10^{-7})$ to $O(10^{-5})$ at $t=300$.
Finally, we computed the relative speed error of the simulations, defined as $(C^{n}-c_s)/c_s$, where $C^{n}=(x^{\ast}(t^{n}+\delta t)-x^{\ast}(t^{n}))/\delta t$  and $x^{\ast}$ an approximation of the center of the pulse, i.~e. the position of its most negative excursion. When we choose $\delta t=1$ the absolute values of the specific error never exceeded $5\times 10^{-15}$ for all $\gamma$; the mean value of the speed remained constant during the computations.

\begin{figure}[!htbp]
\centering
{\includegraphics[width=\textwidth]{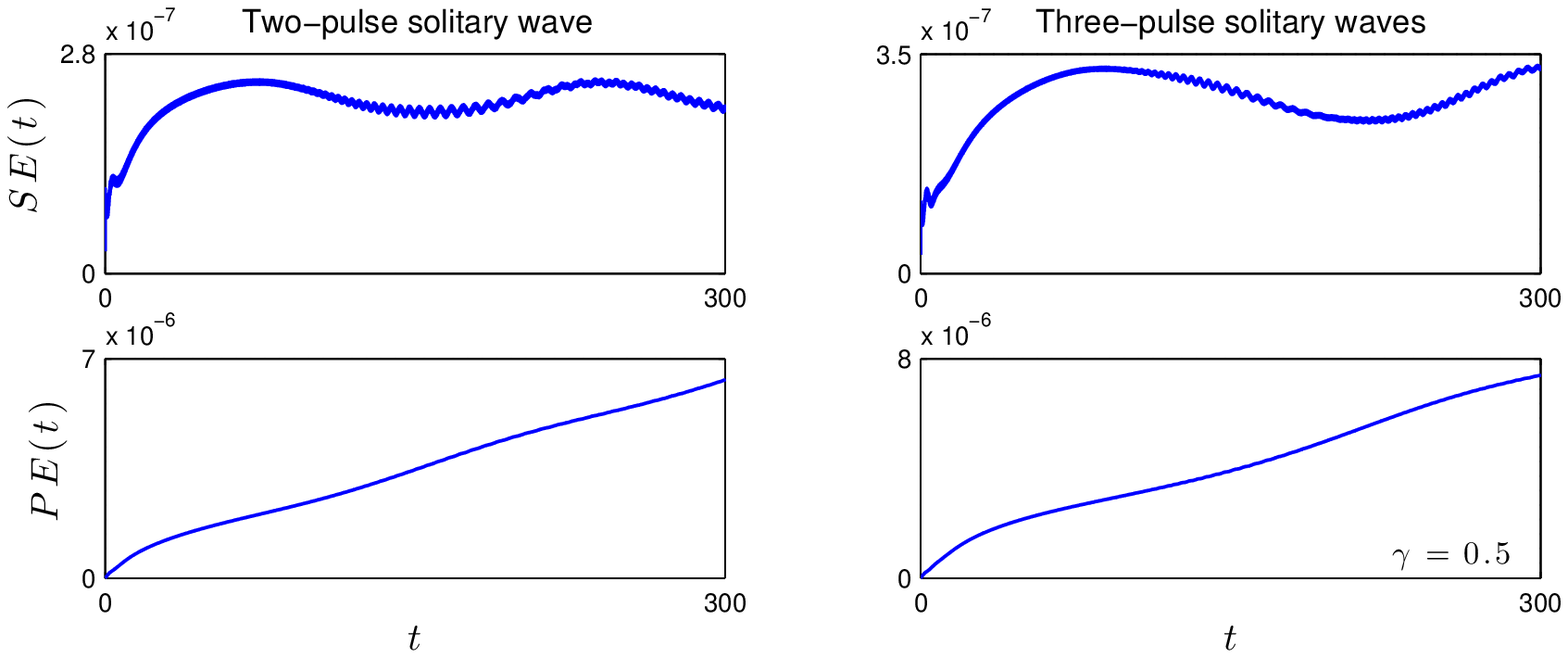}}
\caption{Shape and phase error of the numerical propagation of multi-pulse solitary waves with $c_s=0.75$, $\gamma=0.5$.}%
\label{F10}
\end{figure}

Finally, as a measure of the quality of the numerically generated travelling multi-pulse solitary waves, we present in Figures \ref{F10}, the shape and phase errors during the numerical propagation of two--pulse and three--pulse solitary waves with $c_{s}=0.75$ and $\gamma=0.5$. The shape errors are of $O(10^{-7})$ while the phase errors of about $O(10^{-5})$ at $t=300$.

In conclusion, the outcome of the numerous tests performed in this and the preceding section of the validity and accuracy of the numerical technique for generating initial solitary-wave profiles and of the fully discrete hybrid scheme that was used for their numerical evolution, give us enough confidence to use these schemes in the study of interactions and stability of solitary waves of the Benjamin equation to be undertaken presently.

\section{Overtaking collisions of solitary waves}
\label{sec4}
In this section we study in some detail, by computational means and using the hybrid method, {\it overtaking collisions} of solitary waves of the Benjamin equation. For a given value of $\gamma\in (0,1)$ solitary waves with smaller (absolute) amplitude (i.~e. a smaller in absolute value maximum negative excursion) have larger speed and will consequently overtake solitary waves with larger (absolute) amplitude, which are slower. The solitary waves interact nonlinearly and emerge largely unchanged; their interaction is inelastic, i.~e. it is accompanied by the production of a small amplitude {\it dispersive} tail since the Benjamin equation does not appear to be completely integrable, as already noted in \cite{KB} where results of a simulation of an overtaking collision for solitary waves of the Benjamin equation have been shown.

To set the stage we first present, as a benchmark, the results of a simulation with the hybrid method of an overtaking collision of two solitary waves of the BO equation. The initial solitary waves (cf. (\ref{E317})) had amplitudes $A_{1}=4, A_{2}=1$ and corresponding speeds $c_{s,1}=2$ and $c_{s,2}=1.25$ and were centered at $x_{0,1}=-100$ and $x_{0,2}=100$, respectively. The computation was effected with $N=4096$ and $k=h/20$ on $[-256,256]$, and produced the evolution depicted in Figures \ref{F12}--\ref{F13} at selected instances of $t\in [0,400]$. The two solitary waves interact elastically around $t=265$. During the interactions there always are two distinct peaks present. No artificial oscillations accompany the numerical solution after the interaction
\begin{figure}[!htbp]
\centering
{\includegraphics[width=\textwidth]{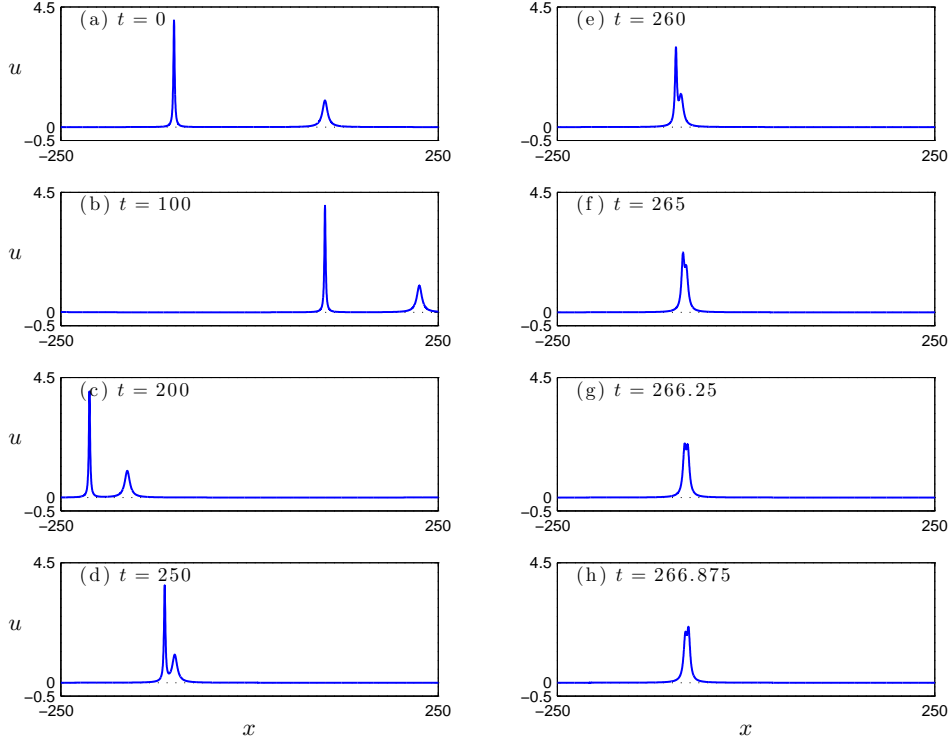}}
\caption{Overtaking collision of two solitary waves of the Benjamin-Ono equation.}%
\label{F12}
\end{figure}

\begin{figure}[!htbp]
\centering
{\includegraphics[width=\textwidth]{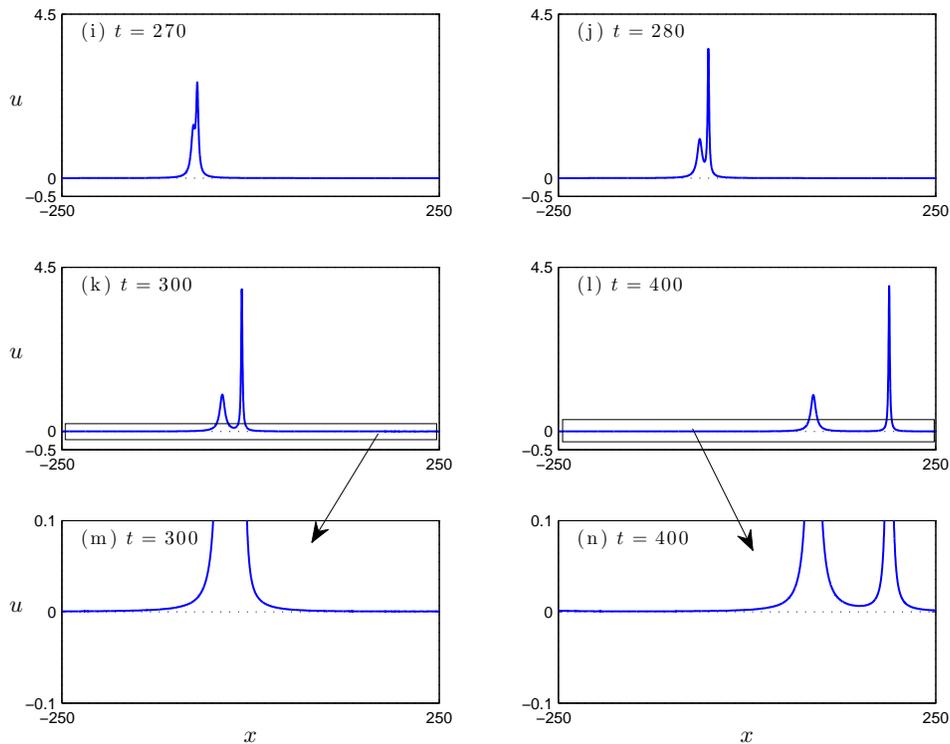}}
\caption{Continuation of results in Figure \ref{F12}. (The profiles (m) and (n) are magnifications of (k) and (l), respectively.)}%
\label{F13}%
\end{figure}

We now turn to the simulations of overtaking collisions of pairs of solitary waves of the Benjamin equation. We studied such collisions for various values of $\gamma$; we present here the results for $\gamma=0.1$ and $\gamma=0.99$. For all cases we used the hybrid method on the spatial interval $[-512,512]$ with $h=0.125$ and $k=0.02$ and constructed initial solitary-wave profiles of various speeds (centered at $x_{1}=256$ and $x_{2}=-256$) by the procedure described in Section \ref{sec3}.

\begin{figure}[!htbp]
\centering
{\includegraphics[width=\textwidth]{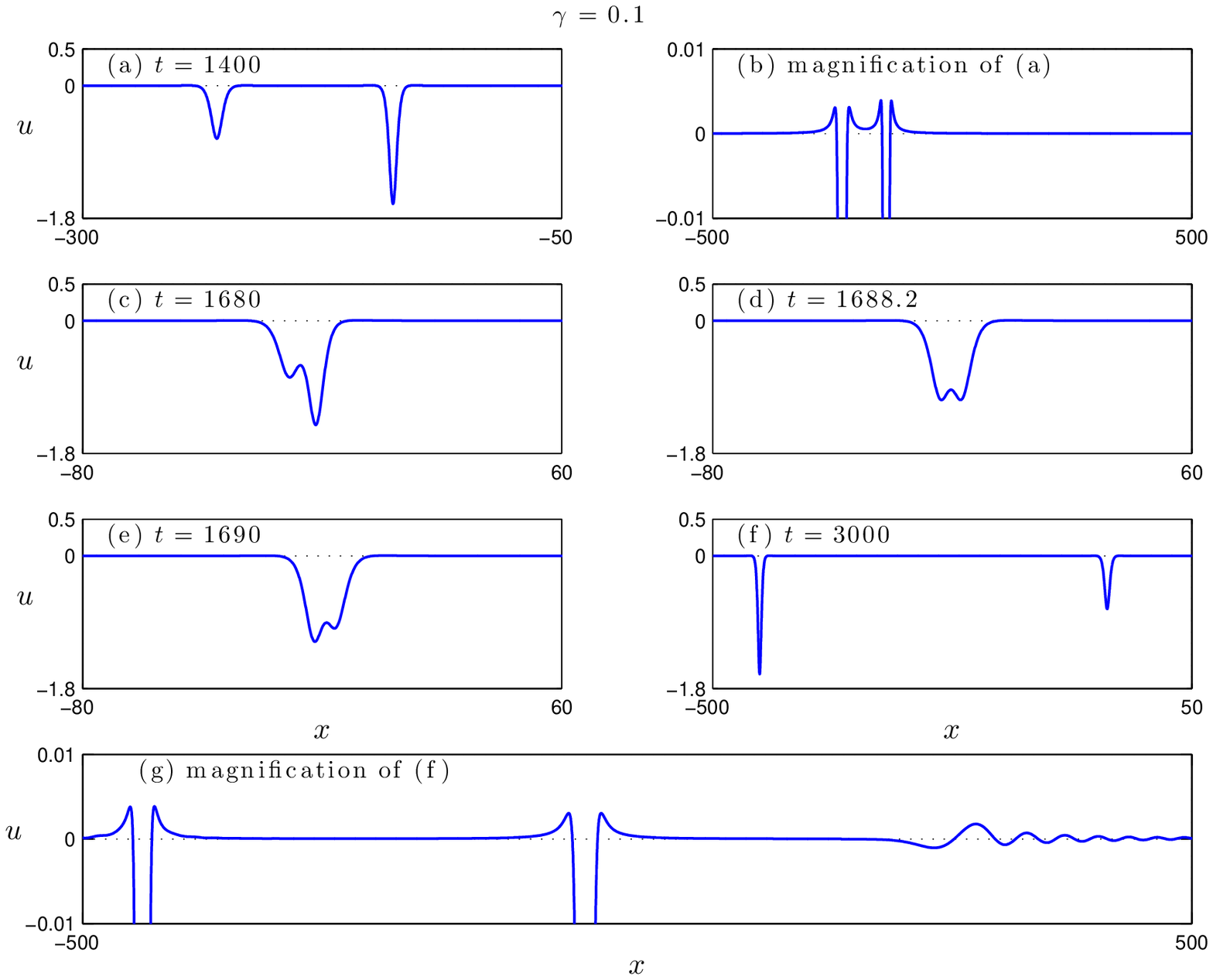}}
\caption{Overtaking collision of solitary waves of the Benjamin equation for $\gamma=0.1$, $c_{s,1}=0.45$, $c_{s,2}=0.75$.}%
\label{F14}%
\end{figure}

Figure \ref{F14} shows several temporal instances of the overtaking collision of two solitary waves of speeds $c_{s,1}=0.45$ and $c_{s,2}=0.75$ in the case $\gamma=0.1$. (During this simulation the $L^{2}$ norm of the solution was $||u||=3.387194802$, and the value of the invariant quantity $H=I+E$ was $H=4.04751039$ up to $T=3000$.) The faster solitary wave overtakes the slower and they interact nonlinearly with two peaks always present during the interaction. The collision produces a dispersive tail (see Figure \ref{F14}(g)), a fact suggesting that the Benjamin equation is not integrable.
\begin{figure}[!htbp]
\centering
{\includegraphics[width=\textwidth]{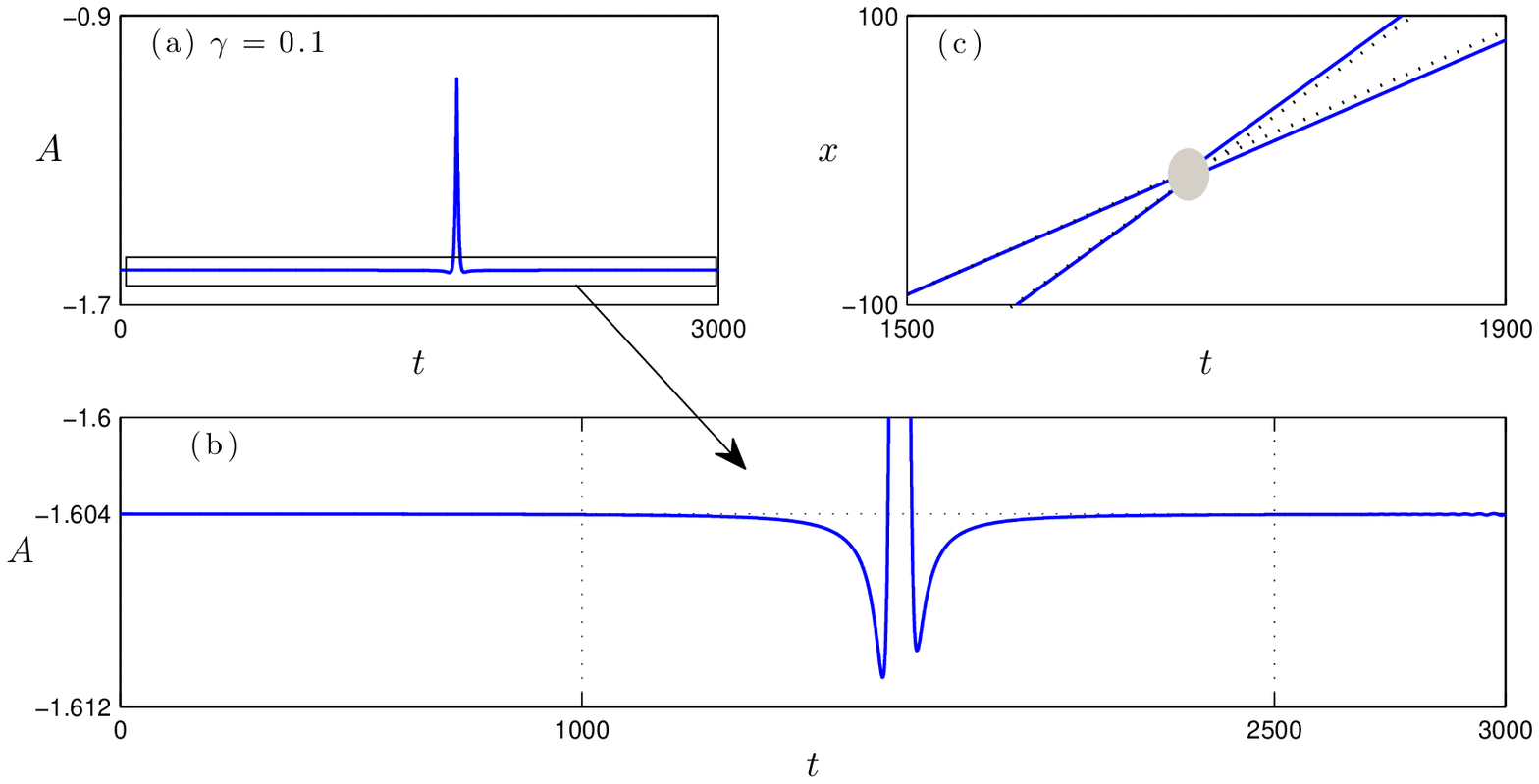}}
\caption{Overtaking collision of solitary waves of the Benjamin equation for $\gamma=0.1$, $c_{s,1}=0.45$, $c_{s,2}=0.75$. Evolution of Figure \ref{F14}. (a): Temporal evolution of the maximum negative excursion of the solution. (c): Paths of solitary waves. The dotted lines would be the paths if no interactions occurred.}%
\label{F15}%
\end{figure}
Note that the dispersive tail precedes the solitary waves being of smaller amplitude and hence faster in our framework. Figure \ref{F15} shows some details of the interaction: In (a) the maximum negative excursion of the solution is plotted versus time. In (b)--a magnification of (a)--one may observe how the maximum negative excursion of the faster wave approaches asymptotically its initial value. The paths of the solitary waves are plotted in (c): The faster wave is shifted slightly forward and the slower backward after the interaction.

\begin{figure}[!htbp]
\centering
{\includegraphics[width=\textwidth]{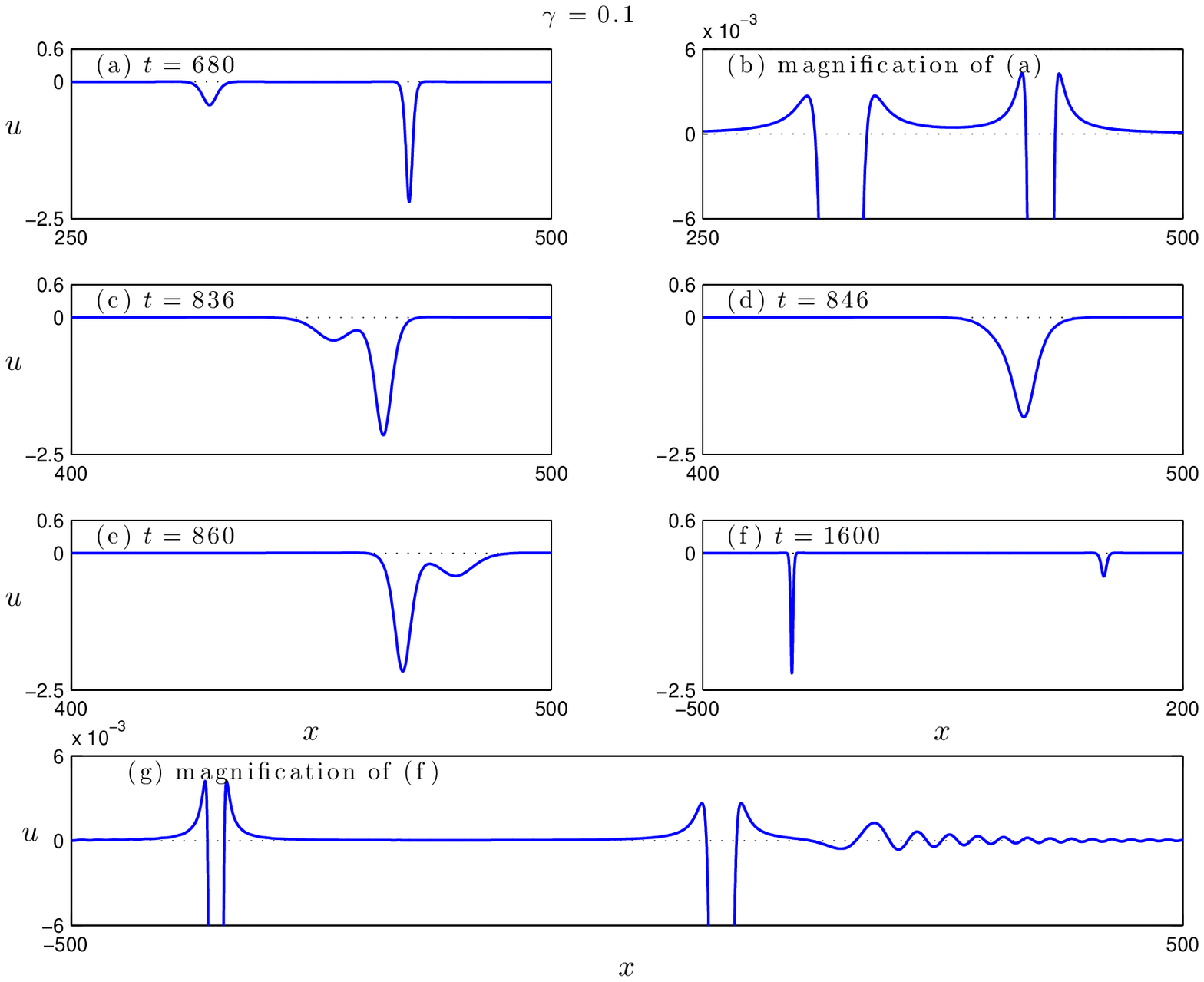}}
\caption{Overtaking collision of solitary waves of the Benjamin equation for $\gamma=0.1$, $c_{s,1}=0.25$, $c_{s,2}=0.85$.}%
\label{F16}%
\end{figure}

\begin{figure}[!htbp]
\centering
{\includegraphics[width=\textwidth]{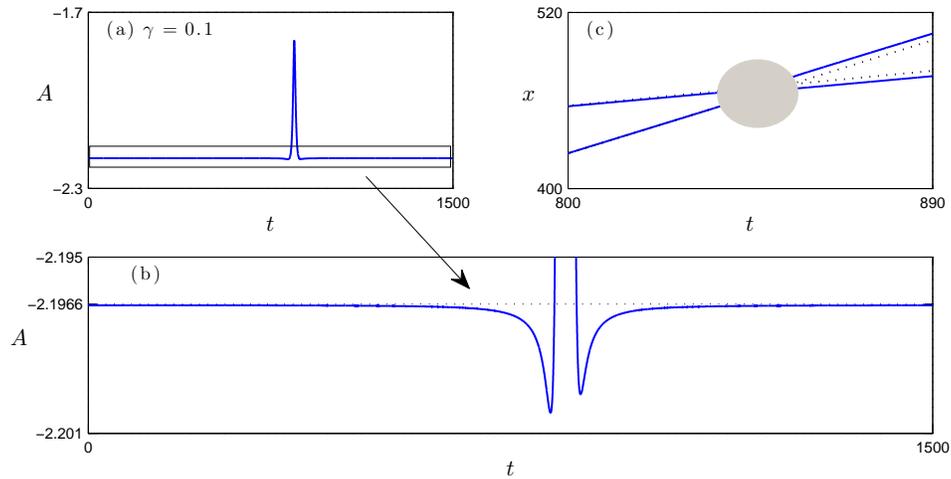}}
\caption{Overtaking collision of solitary waves of the Benjamin equation for $\gamma=0.1$, $c_{s,1}=0.25$, $c_{s,2}=0.85$. Graphs analogous to those of Figures \ref{F15}.}%
\label{F17}%
\end{figure}

In Figures \ref{F16}-\ref{F17} we show the analogous simulation of the overtaking collision of two solitary waves of the Benjamin equation of initial speeds $c_{s,1}=0.25$ and $c_{s,2}=0.85$, again for $\gamma=0.1$. The larger difference of the speeds in this experiment apparently causes the formation of a single peak momentarily during the interaction. Otherwise the details of the overtaking collision are qualitatively the same with those in Figures \ref{F14}-\ref{F15}. During this simulation the values of the invariants $||u||$ and $H$ remained equal to $3.93689569$ and $4.42223526$, respectively, up to $T=1500$.

\begin{figure}[!htbp]
\centering
{\includegraphics[width=\textwidth]{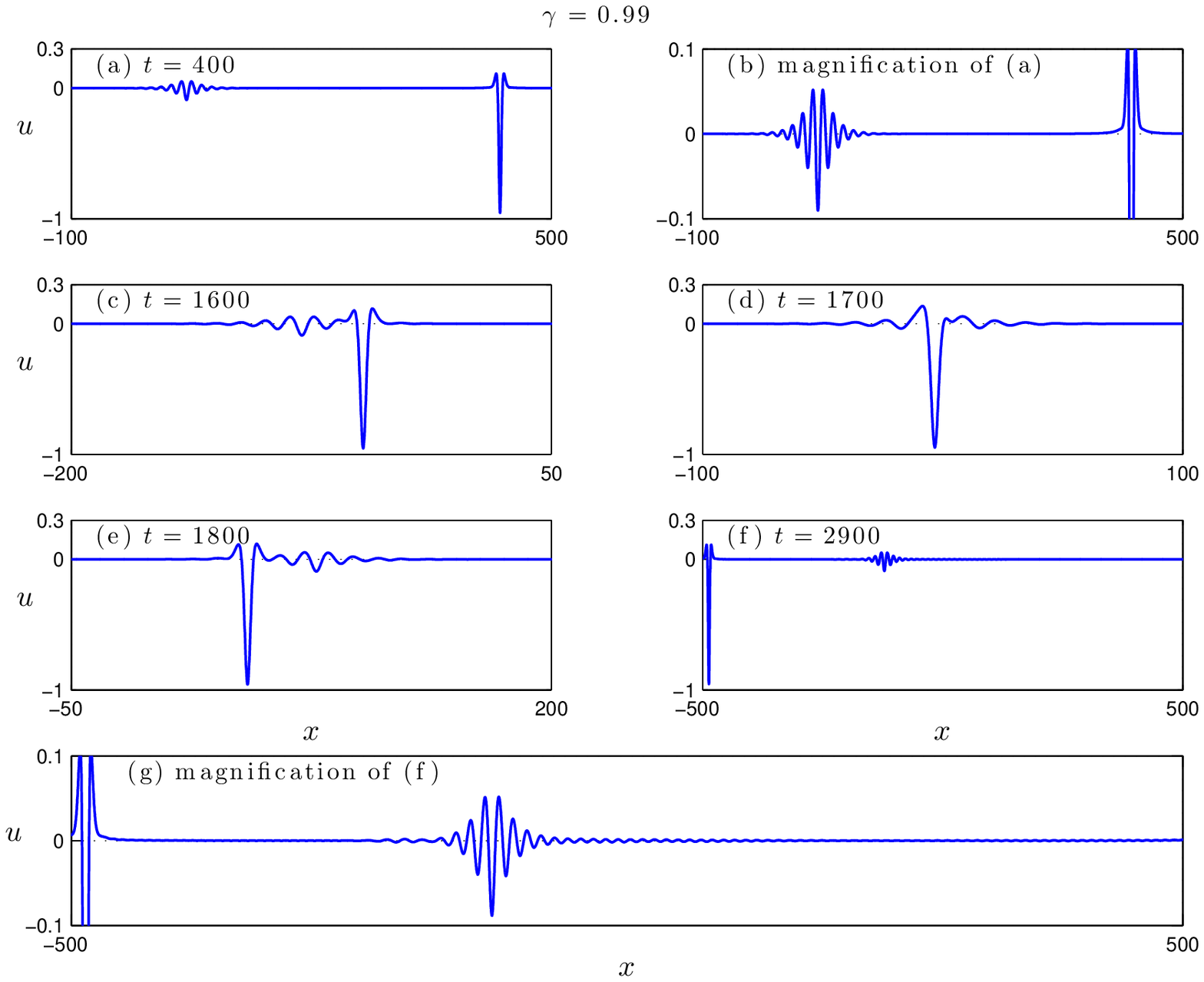}}
\caption{Overtaking collision of solitary waves of the Benjamin equation for $\gamma=0.99$, $c_{s,1}=0.45$, $c_{s,2}=0.75$.}%
\label{F18}%
\end{figure}

We noticed that the collisions became harder to simulate for $\gamma>0.9$. Figure \ref{F18} shows the interaction of two solitary waves of speeds $c_{s,1}=0.45$ and $c_{s,2}=0.75$ in the case $\gamma=0.99$. The $L^{2}$ norm was preserved to ten digits (it was equal to $1.532051456$) up to $t=3000$, but $H=6.821038$ was preserved to $7$ digits, reflecting the increased difficulty of the computation. It is not clear whether the small oscillations in front of the smaller, highly oscillatory solitary wave in Figure \ref{F18}(g) at $t=2900$ belong to a dispersive tail or are numerical artifacts or somehow indicate that the smaller wave has not yet stabilized after the interaction.
\begin{figure}[!htbp]
\centering
{\includegraphics[width=\textwidth]{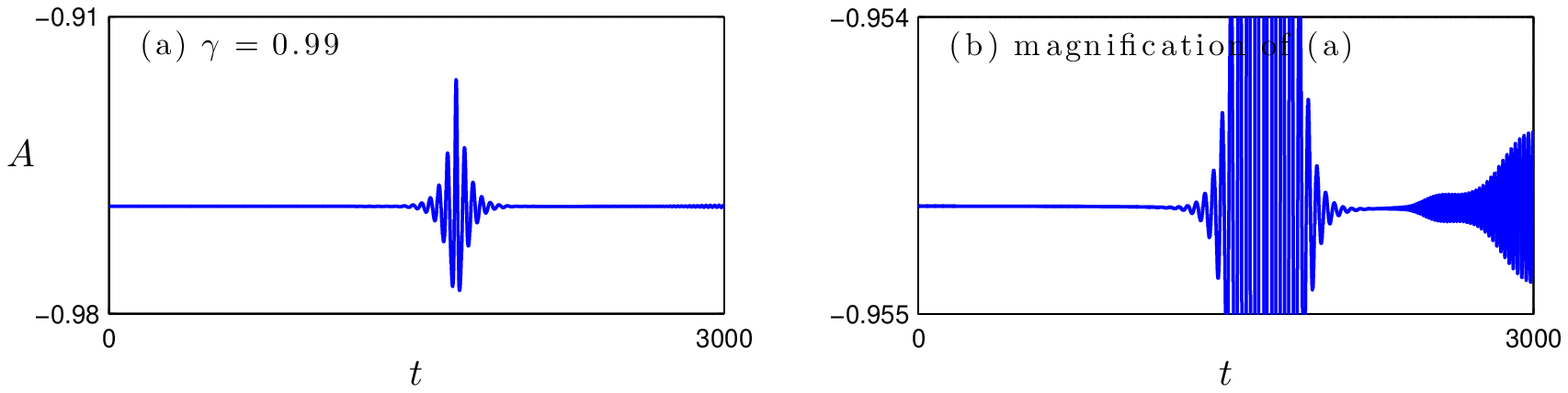}}
\caption{Overtaking collision of solitary waves of the Benjamin equation for $\gamma=0.99$, $c_{s,1}=0.45$, $c_{s,2}=0.75$. Graphs analogous to (a) and (b) of Figure \ref{F15}.}%
\label{F19}%
\end{figure}
We observe that after about $t=2500$ as shown in Figure \ref{F19} in which the maximum negative excursion of the solution is plotted versus time, after achieving again its pre-interaction value, the maximum negative excursion of the slower wave starts oscillating as it interacts with the dispersive tail.

\begin{figure}[!htbp]
\centering
{\includegraphics[width=\textwidth]{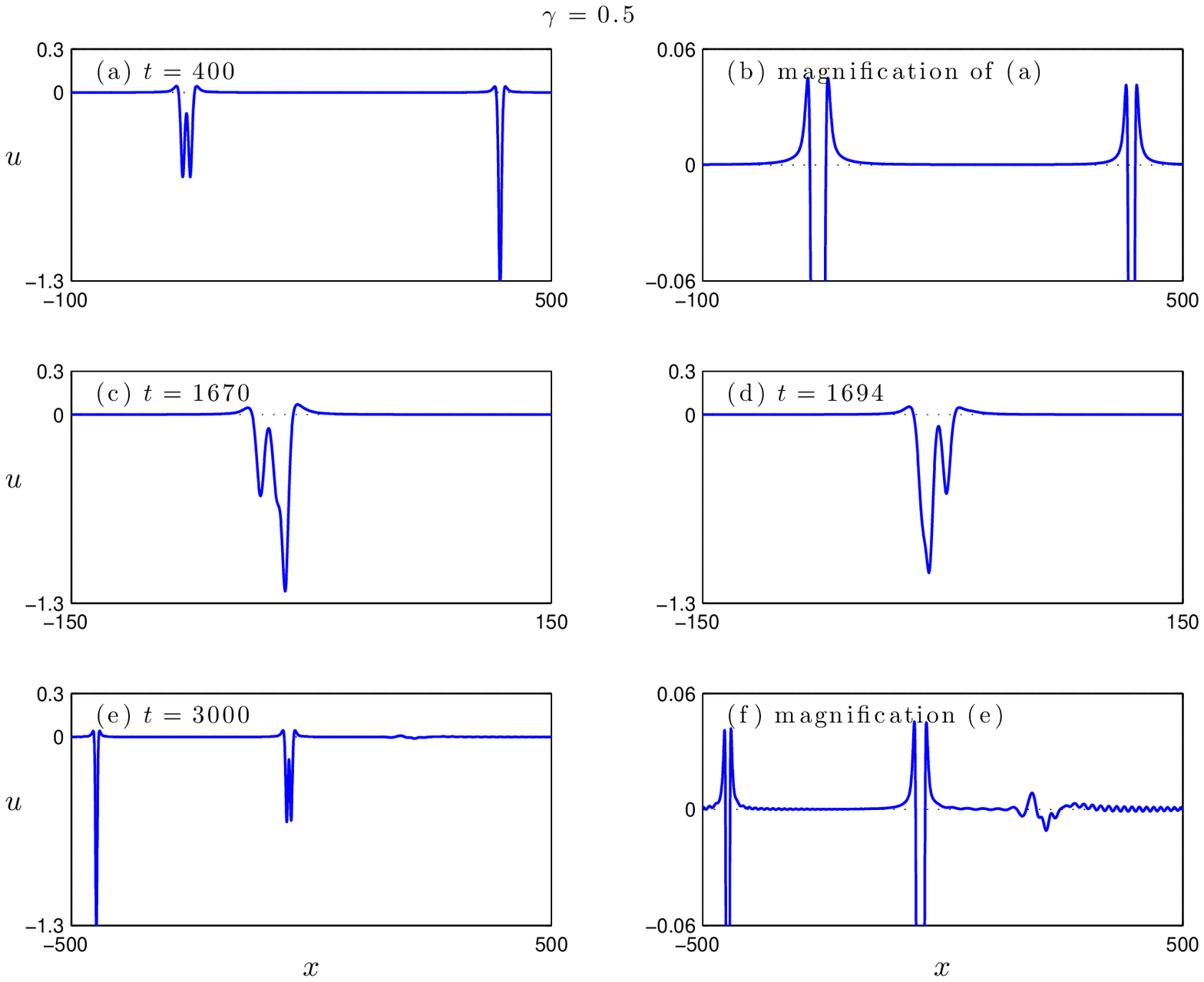}}
\caption{Overtaking collision of a two-pulse and an ordinary solitary wave of the Benjamin equation for $\gamma=0.5$, $c_{s,1}=0.45$, $c_{s,2}=0.75$.}%
\label{F20}%
\end{figure}
\begin{figure}[!htbp]
\centering
{\includegraphics[width=\textwidth]{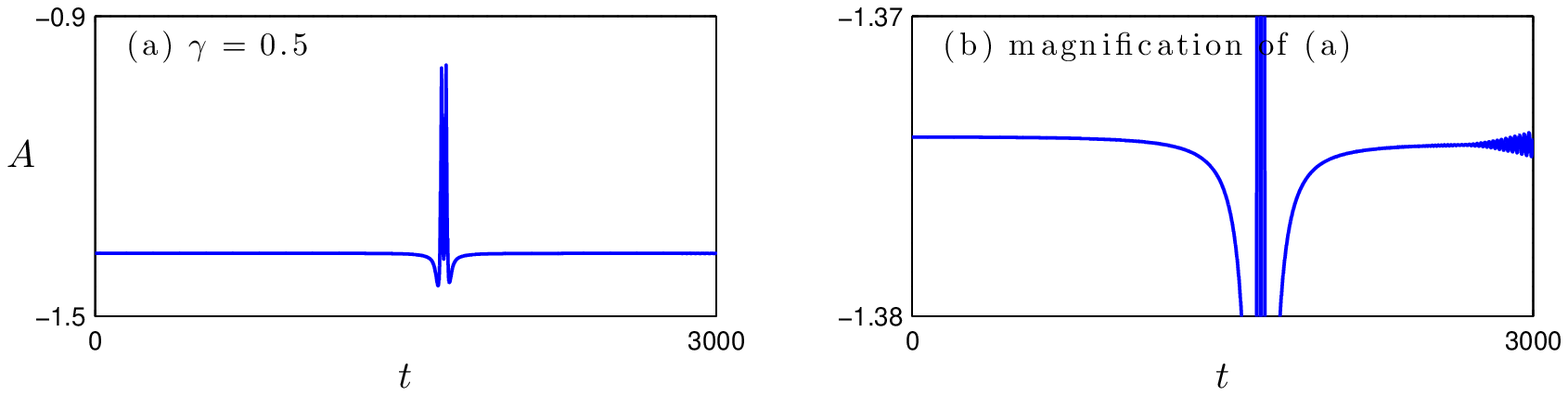}}
\caption{Overtaking collision of solitary waves for the Benjamin equation for $\gamma=0.5$, $c_{s,1}=0.45$, $c_{s,2}=0.75$, evolution of Figure \ref{F20}. Maximum negative excursion of the solution versus time.}%
\label{F21}%
\end{figure}

We also performed numerical experiments simulating overtaking collisions involving multi-pulse solitary waves of the Benjamin equation. Figures \ref{F20} and \ref{F21} show such an interaction of a fast two-pulse solitary wave of speed $c_{s,2}=0.75$ with a slower single-pulse wave with $c_{s,1}=0.45$ for $\gamma=0.5$. During this simulation we observed that $||u||=2.873492446, H=2.8836586$ up to $t=3000$. After the interaction the waves separate and there is evidence of a dispersive tail, but the two-pulse wave has not quite recovered its shape and initial amplitudes by $t=3000$. The same is true for the single-pulse wave whose maximum negative excursion has not returned to its initial value by $t=3000$ as Figure \ref{F21} indicates.

\section{Stability of solitary waves}
\label{sec5}
In this section we first study by computational means the stability of single- and multi-pulse solitary waves of the Benjamin equation under small perturbations. As was mentioned in the Introduction, a theory of stability of single-pulse waves was outlined in \cite{B2} and a complete proof for small $\gamma$ was given in \cite{ABR}. Another proof, valid for all $\gamma\in [0,1)$, of stability in a weaker sense was given in \cite{A}.

\begin{figure}[!htbp]
\centering
{\includegraphics[width=\textwidth]{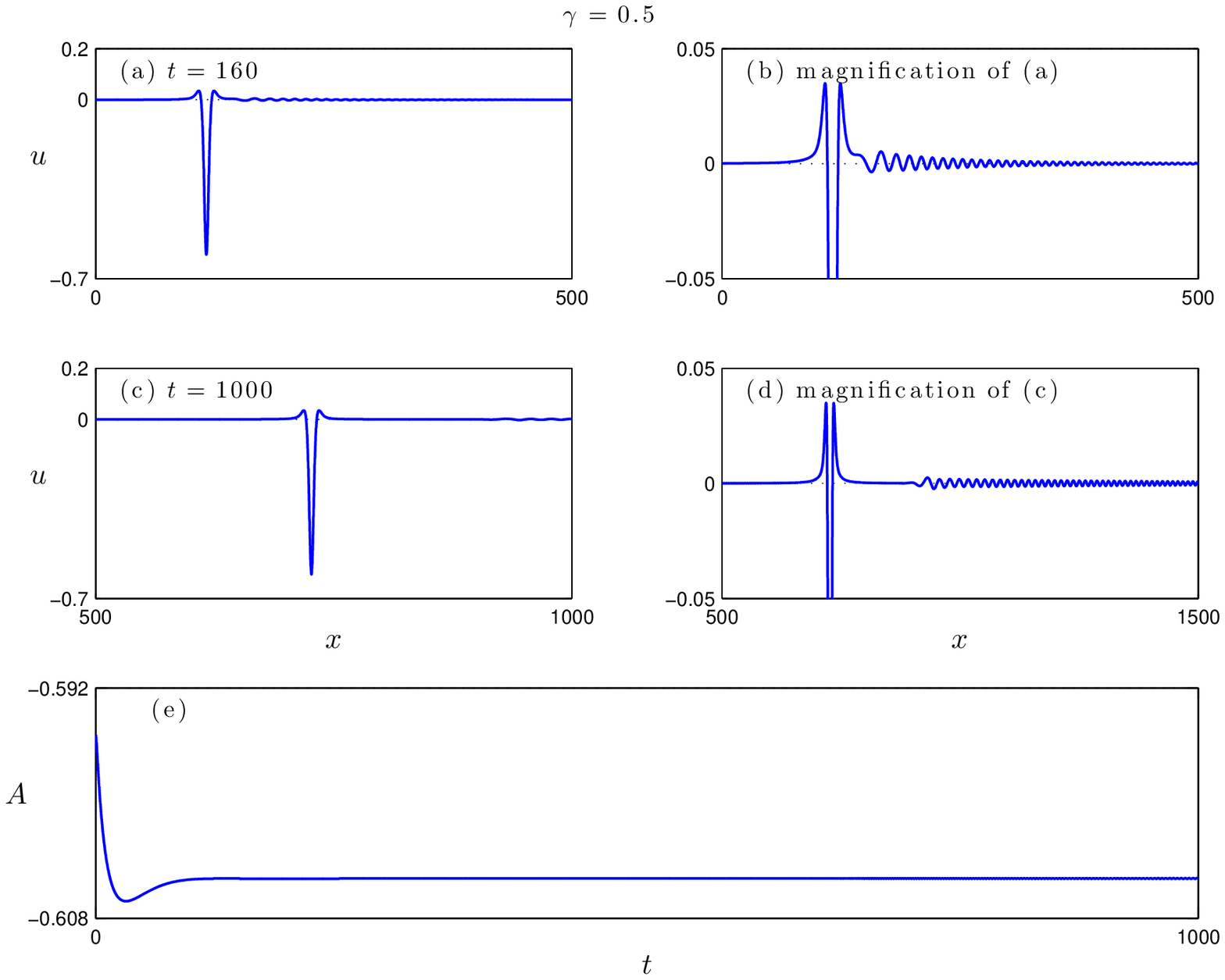}}
\caption{Evolution of a perturbed single-pulse solitary wave of the Benjamin equation ($\gamma=0.5$). (b) and (d) are magnifications of (a) and (c), respectively. (e): Evolution of the maximum negative excursion of the solution.}%
\label{F22}%
\end{figure}

\begin{figure}[!htbp]
\centering
{\includegraphics[width=\textwidth]{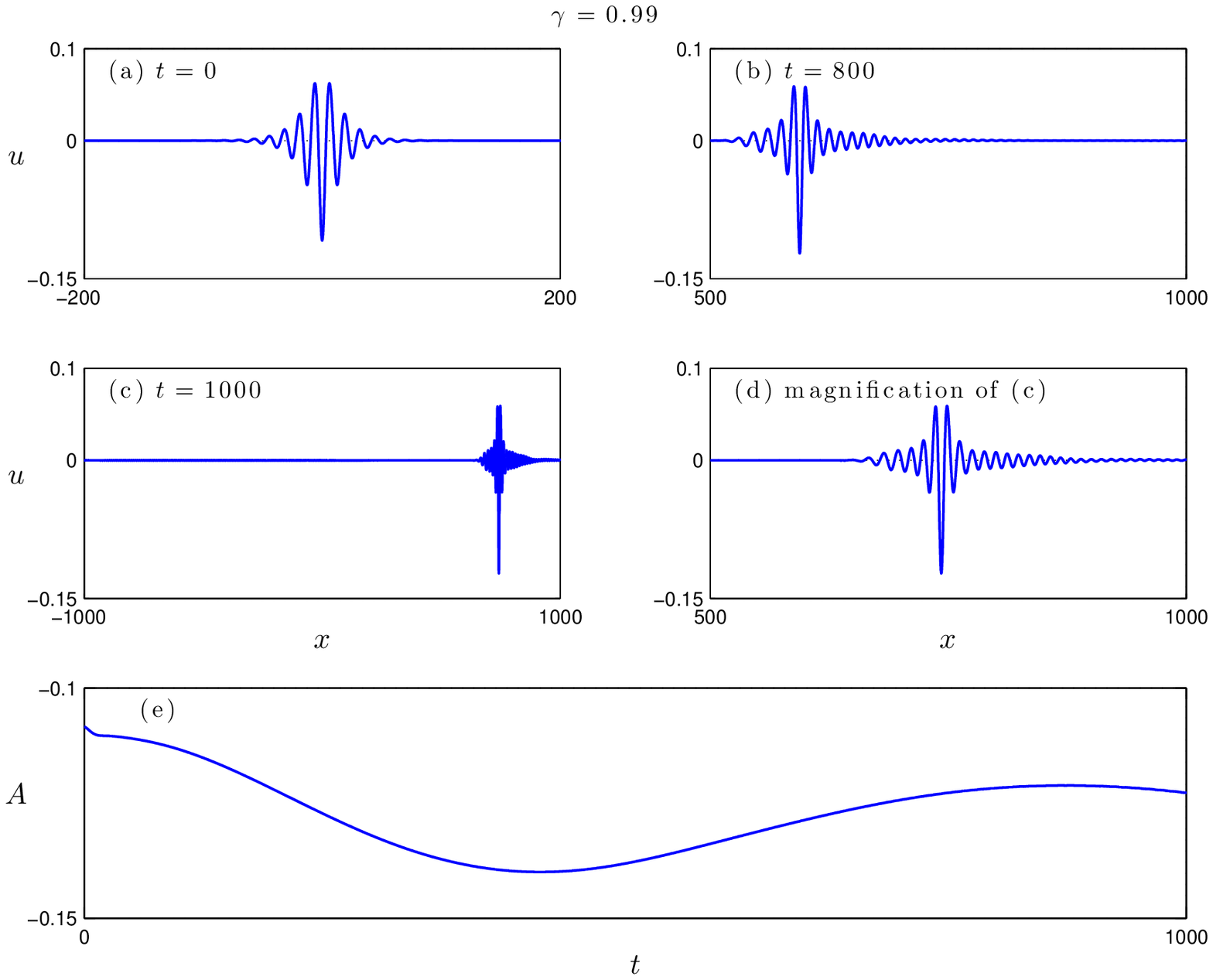}}
\caption{Evolution of a perturbed solitary wave of the Benjamin equation ($\gamma=0.99$). ((b) and (d) are magnifications of (a) and (c), respectively.) (e): Evolution of the maximum negative excursion of the solution.}%
\label{F23}%
\end{figure}

We start with the single-pulse case. Figure \ref{F22}(a)--(d) shows the evolution (effected with the hybrid method on the spatial interval $[-2048,2048]$ with $h=0.0625$ and $k=0.02$) ensuing from a single-pulse solitary wave with $\gamma=0.5$ and $c_{s}=0.75$, centered at $x_{0}=0$, when it is perturbed by a multiplicative factor $r=1.1$. As expected, the perturbed solitary wave evolves into a new one of slightly larger maximum negative excursion plus a preceding dispersive tail. Figure \ref{F22}(e) shows the evolution of the maximum negative excursion of the solitary wave from its initial value $-0.59526$ to its eventual value which is equal to $-0.60523$.
We also simulated the evolution of a perturbed solitary wave corresponding to $\gamma=0.99$. Figure \ref{F23}(a)--(d) shows this evolution. The initial solitary wave had $c_{s}=0.75$ and was perturbed by a multiplicative factor of $r=1.2$. (The computation was effected on $[-1024,1024]$ with $h=0.0625, k=0.02$ up to $T=1000$.) The wave radiates forward a small-amplitude oscillatory wavetrain which has not separated from the main wave up to $T=1000$. This fact, and also the temporal variation of the maximum negative excursion of the wave (Figure \ref{F23}(e)) which has not achieved an asymptotic state by $t=1000$, does not allow us to reach a conclusion about the stability of solitary waves for $\gamma=0.99$. The wave may be unstable and keep radiating small-amplitude oscillations for all $t$ or may stabilize into a nearby solitary wave after very long time.

\begin{figure}[!htbp]
\centering
{\includegraphics[width=\textwidth]{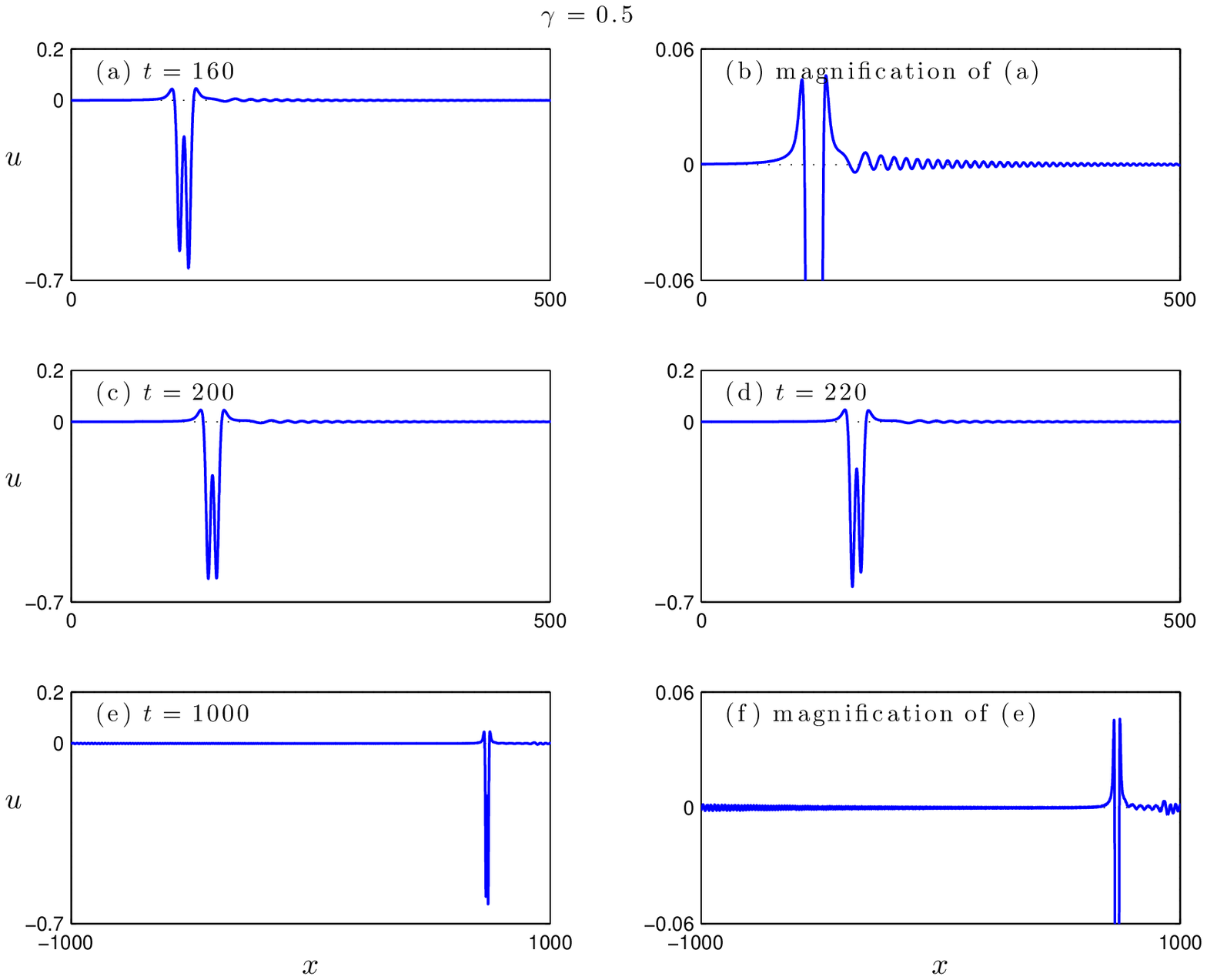}}
\caption{Evolution of a perturbed two-pulse solitary wave of the Benjamin equation. ($\gamma=0.5$. ((b) and (f) are magnifications of (a) and (e), respectively.))}%
\label{F24}%
\end{figure}

\begin{figure}[!htbp]
\centering
{\includegraphics[width=\textwidth]{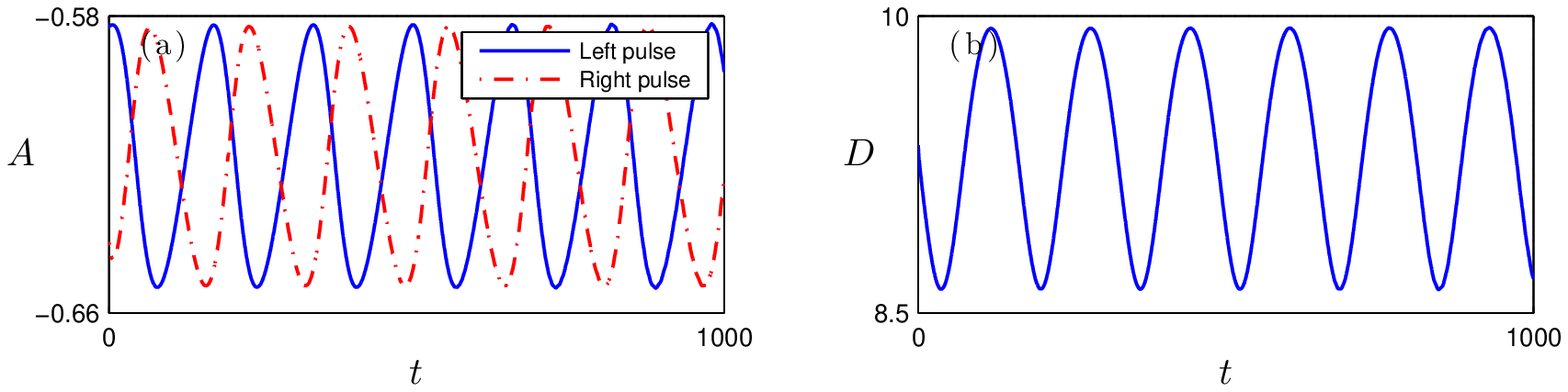}}
\caption{(a): Amplitudes (maximum negative excursions) of the two negative peaks of the perturbed two-pulse solitary wave of Figure \ref{F24}, and (b): Distance between the two peaks, as functions of $t$.}%
\label{F25}%
\end{figure}

We turn now to a stability study of a two-pulse solitary wave. We took as initial condition a two-pulse solitary wave in the case $\gamma=0.5$ and perturbed it asymmetrically multiplying it by a factor $r({\rm tanh}x +1)+1$ with $r=0.05$. Figure \ref{F24} shows the evolution that ensues. (The computation was done on $[-1024,1024]$ up to $T=1000$ using $h=0.0625, k=0.02$.) The perturbed two-pulse wave radiates forward the usual small-amplitude oscillatory wavetrain. We observe that its two negative peaks oscillate exchanging heights in a periodic-like manner (Figure \ref{F25}(a)), while their distance is also oscillating apparently periodically (Figure \ref{F25}(b)). This \lq dance\rq\ of the twin peaks went on up to the end of our computation at $t=1000$, but it is unlikely to continue unaltered for ever due to the constant shedding of radiation.

\begin{figure}[!htbp]
\centering
{\includegraphics[width=\textwidth]{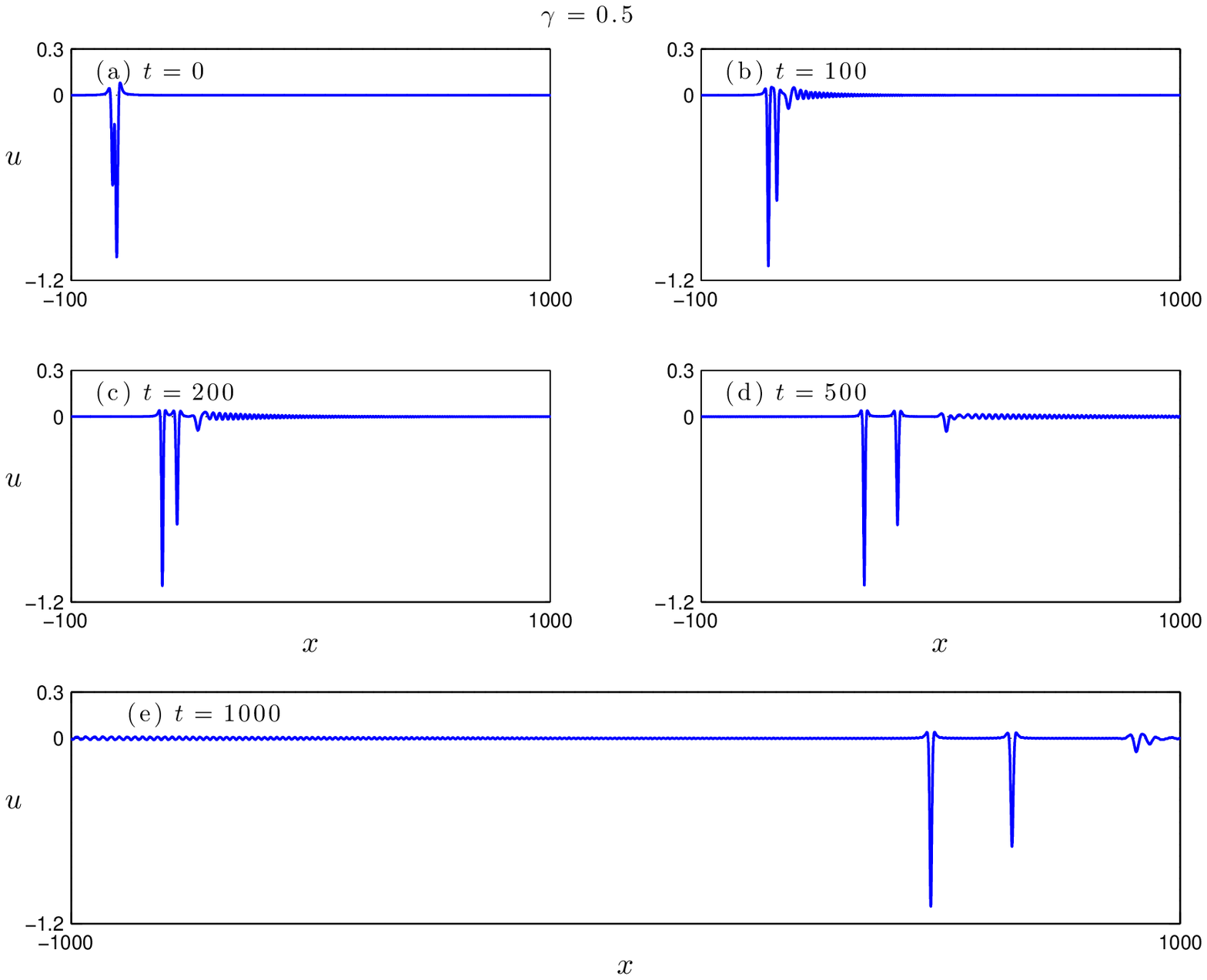}}
\caption{Evolution of a more perturbed two-pulse solitary wave of the Benjamin equation ($\gamma=0.5$). }
\label{F26}%
\end{figure}

\begin{figure}[!htbp]
\centering
{\includegraphics[width=\textwidth]{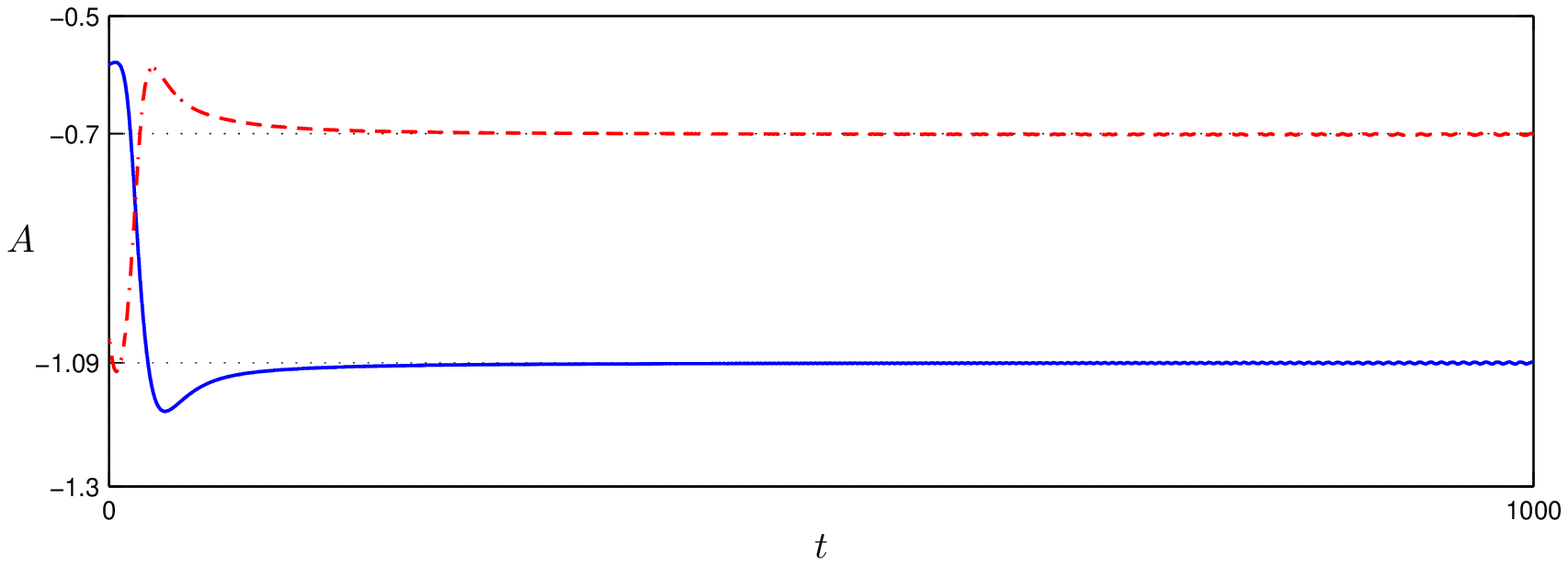}}
\caption{Maximum negative excursions of the two negative peaks as functions of $t$.}\label{F27}
\end{figure}

In a related numerical experiment, whose outcome is shown in Figure \ref{F26}, we perturbed the same initial two-pulse solitary wave with a larger asymmetric factor ($r$ was taken now to be $0.4$) of the same form as above. (All computational parameters remained the same.)
After a brief initial dancing phase (up to about $t=40$) accompanied by radiation, we observed that two single-pulse solitary waves were generated. Figure \ref{F27} shows the evolution of the maximum negative excursions of the two negative peaks up to $T=1000$.

We conclude then that the effect of the larger perturbation is apparently to accelerate the end of the dance and initiate resolution into solitary waves.

As was already mentioned in the Introduction, Kalisch and Bona in \cite{KB} describe numerical experiments in which they observed resolution into solitary waves for the Benjamin equation with initial Gaussian profiles of the form $Ae^{-(x/\lambda)^{2}}$. As $\lambda$ was increased the emergence of a pair of \lq orbiting\rq\ solitary waves was observed which danced in the way previously described. For larger values of $\lambda$, they report that \lq triplets\rq\ and \lq quadruplets\rq\ of such solitary waves appeared. It was further conjectured in \cite{KB} (on the basis of the observed increase of the distance between the peaks of the orbiting pairs of solitary waves) that the system `may eventually transform into two separately propagating solitary waves'.

In the light of the numerical experiments of the present paper one could interpret the orbiting solitary waves of \cite{KB} as perturbed multi-pulse solitary waves, which, after an intermediate dancing stage, resolve themselves into separate single-pulse solitary waves.

\begin{figure}[!htbp]
\centering
{\includegraphics[width=\textwidth]{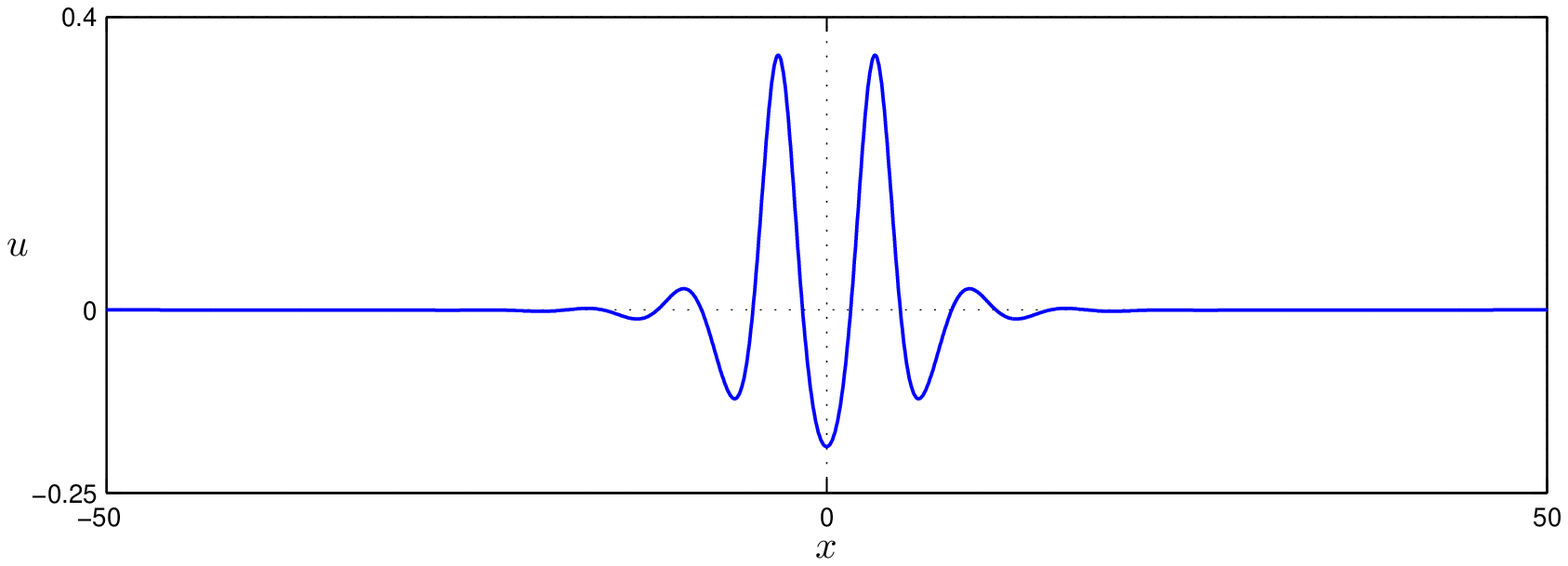}}
\caption{Initial \lq depression\rq\ solitary wave $\phi(x)=u_{0}(x), \gamma=0.94, c_{s}=0.9$.}\label{F28}
\end{figure}

As was mentioned in the Introduction we also computed the evolution of \lq depression\rq\ solitary waves of the Benjamin equation considered in \cite{CA} with the aim of studying their stability properties. In order to facilitate comparisons with the results of \cite{CA}, we computed the initial \lq depression\rq\ wave profile by solving the solitary-wave equation in the form given by equation (44) of \cite{CA}, i.~e. as solution $\phi=\phi(x)$ of
\begin{eqnarray*}
\nu \phi-\phi^{2}-2\gamma \mathcal{H}\phi_{x}-\phi_{xx}=0,
\end{eqnarray*}
with $\nu=1, \gamma=0.94$. For this purpose we used the CGN algorithm (without continuation) taking as initial guess the usual (\lq elevation\rq) solitary wave of the Benjamin equation corresponding to $\gamma=0.94, c_{s}=0.9$, reflected about the $x-$axis and multiplied by a factor of two. (We performed $175$ iterations with a final residual error of the order of $10^{-13}$.) The profile $\phi(x)=u_{0}(x)$ that was obtained is shown in Figure \ref{F28}; it corresponds to the profile of the uppermost snapshot of Figure 6 of \cite{CA}.

\begin{figure}[!htbp]
\centering
{\includegraphics[width=\textwidth]{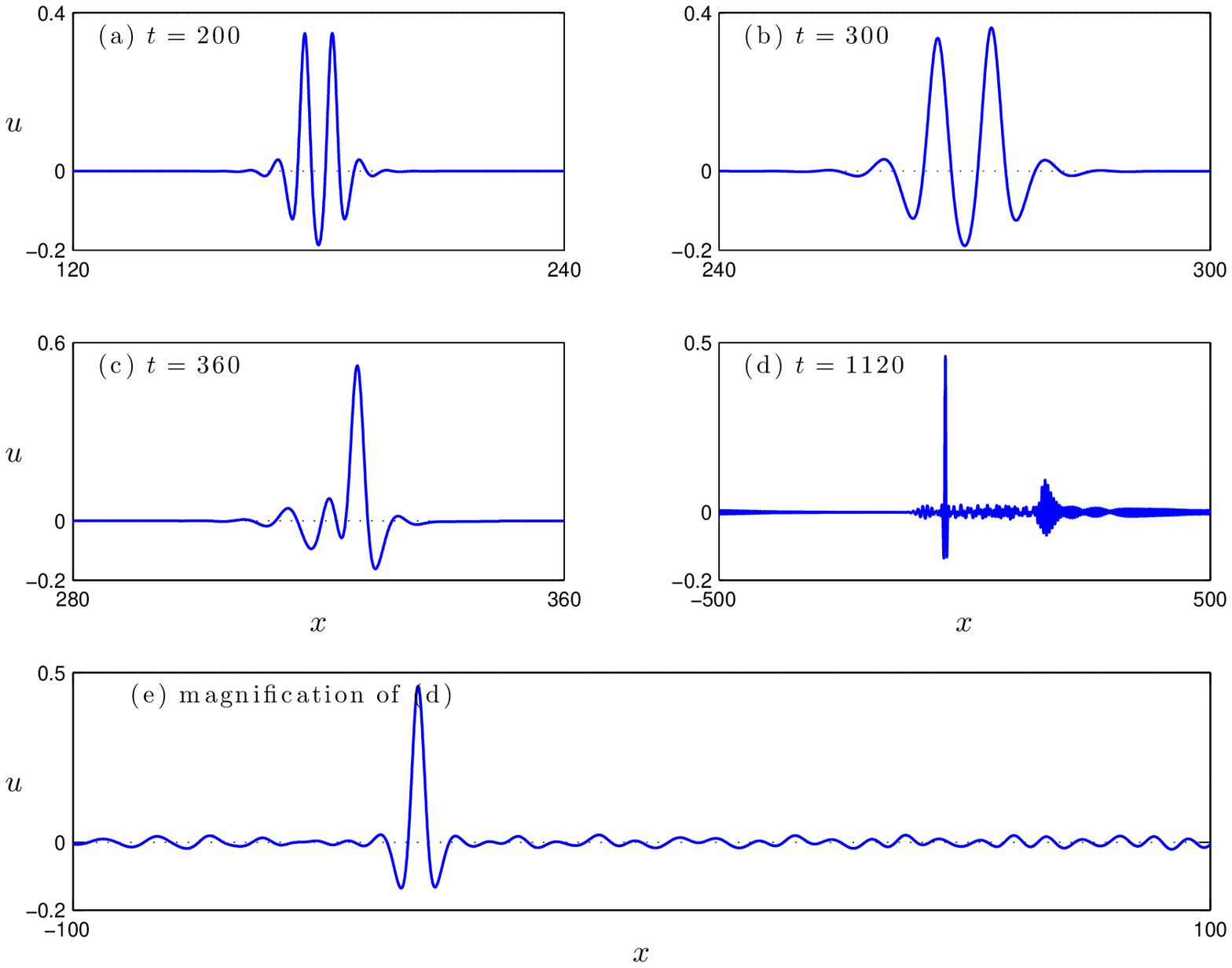}}
\caption{Evolution of the initial \lq depression\rq\ solitary wave profile $u_{0}(x)$ of Figure \ref{F28} under the p.d.e. (\ref{E512}); (e) is a magnification of (d) in the neighborhood of the main pulse.}%
\label{F29}%
\end{figure}

We then integrated forward in time with our hybrid scheme using the appropriate transformed version of the p.d.e. (43) of \cite{CA}. Specifically, if $\eta=\eta(X,\tau)$ is the solution of that equation, our change of variables was defined by
\begin{equation}\label{E511}
\eta(X,\tau)=u(x,t),\quad x=X+2.8\tau, \quad t=2\tau.
\end{equation}
This gave for the variable $u(x,t)$ the Benjamin equation of the form
\begin{equation}\label{E512}
u_{t}+1.4u_{x}-uu_{x}-0.94\mathcal{H}u_{xx}-0.5u_{xxx}=0,
\end{equation}
i.~e. of the form (\ref{E11}) with $\beta=-1, \alpha, \gamma, \delta$ positive, which we integrated with the hybrid method on $[-1024,1024]$ using $h=0.125 (N=16384), k=0.02$ up to $t=1120$. The ensuing evolution is depicted in Figure \ref{F29}.

The initial profile moves to the right with speed $c_{s}=0.9$, apparently unchanged until about $t=250$. (Note that the analogous wave in Figure 6 of \cite{CA} moves to the left because its speed is equal to $-1$. This follows from our change of variables (\ref{E511}) which implies that $u(x,t)=\phi(x-0.9t)$ if and only if $\eta(X,\tau)=\phi(X+\tau)$.) After that time, perturbed by the errors inherent in the numerical scheme the \lq depression\rq\ wave starts losing its shape and eventually develops into one main pulse, apparently a solitary wave of \lq elevation\rq\ , which continues travelling to the right, preceded by a dispersive oscillatory wavetrain. This instability confirms the results of \cite{CA} and may be seen more clearly in another numerical experiment in which we took as initial value the function $ru_{0}(x)$ with $r=1.1$. The evolution that resulted was simulated again up to $t=2200$ with the hybrid scheme for (\ref{E512}) with the same discretization parameters as before and is depicted in Figure \ref{F29}. The perturbed initial \lq depression\rq\ solitary wave loses its shape fast and apparently evolves in two usual (\lq elevation\rq) solitary waves of different heights that travel to the right preceded by a dispersive tail. (Note that in Figures \ref{F28} and \ref{F29} the solitary waves have positive peaks, while in previous sections of the paper at hand they had negative. This is due to the negative sign of the $uu_{x}$ term in (\ref{E512}): If we make the change of variable $v=-u$, $v$ satisfies the Benjamin equation $v_{t}+1.4v_{x}+vv_{x}-0.94\mathcal{H}v_{xx}-0.5v_{xxx}=0$, which is our usual form. For the latter equation the solitary waves of \lq elevation\rq\ type have negative maximum excursions from zero and waves of smaller absolute amplitude are faster than those of larger absolute amplitude, cf. e.~g. Figure \ref{F16}. Hence in the $u-$equation (\ref{E512}) the solitary waves have positive maximum excursions and still move to the right with the waves of smaller amplitude being faster than those of larger amplitude and with the tiny dispersive oscillatory wavetrain being even faster as observed in Figure \ref{F29}.)
\begin{figure}[!htbp]
\centering
{\includegraphics[width=\textwidth]{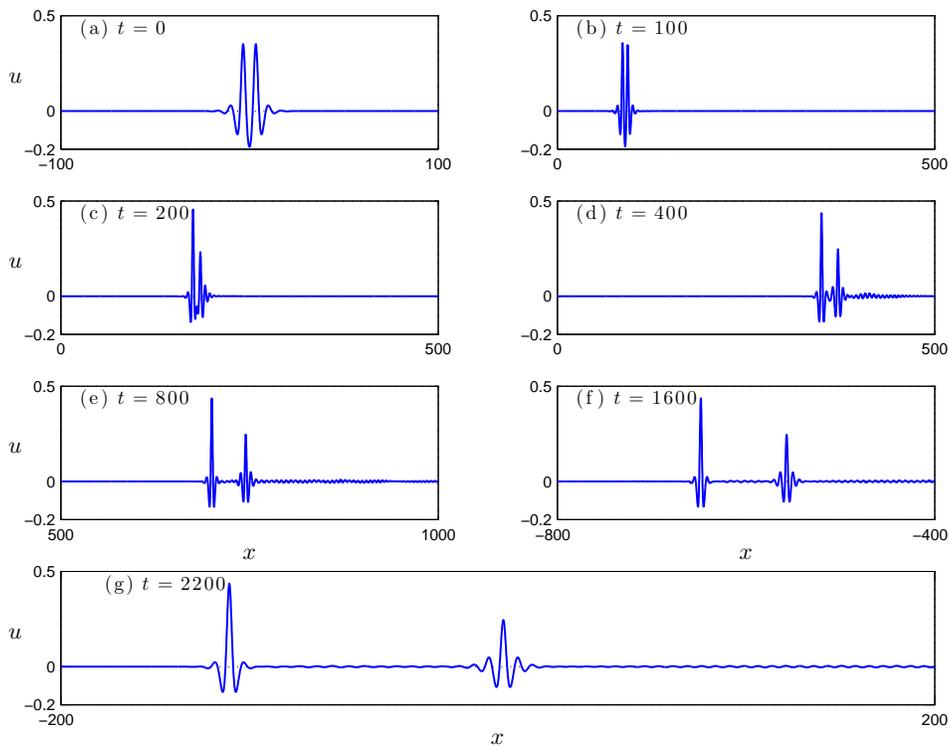}}
\caption{Evolution of the perturbed initial \lq depression\rq\ solitary wave profile $ru_{0}(x)$ (r=1.1, $u_{0}(x)$ as in Figure \ref{F28}) under the p.d.e. (\ref{E512})}%
\label{F30}%
\end{figure}

%%%%%%%%%%%%%%%%%%%%%%%%%%%%%%%%%%%%%%%%%%%%%%%%%%%%%%%%%%

\section*{Acknowledgments}
V. Dougalis and A. Duran have been supported by  project MTM2010-19510/MTM (MCIN).
%%%% Bibliography  %%%%%%%%%%
%\bibliographystyle{alpha}
%\bibliography{biblio}
%\input{bibliog}

\end{document}